\documentclass[11pt]{article}
\usepackage[T1]{fontenc}
\usepackage[utf8]{inputenc}
\usepackage{vien}
\usepackage{mathtools, graphicx, caption, subcaption, bm, epstopdf, url}
\usepackage{algorithm, enumitem}
\usepackage[top=2.54cm,right=3cm,left=3cm,bottom=2.54cm]{geometry}
\usepackage[final]{changes}
\newcommand{\statrv}{S}
\newcommand{\statval}{s}
\newcommand{\batchsize}{m}
\newcommand{\stepsize}{\alpha}
\newcommand{\mmt}{\beta}
\newcommand{\thres}{\gamma}
\newcommand{\subdiff}{\partial}
\newcommand{\opt}{^\star}
\newcommand{\Xopt}{\mc{X}\opt}

\newcommand{\gproj}[2]{{\mathop{\Pi}}_{#1} \left(#2\right)}
\newcommand{\clip}[1]{\mathrm{clip}_{\gamma}\left(#1\right)}
\newcommand{\tvclip}[2]{\mathrm{clip}_{#1}\left(#2\right)}

\title{Stability and Convergence of Stochastic Gradient Clipping:  Beyond Lipschitz Continuity and Smoothness}

\author{	
	Vien V.~Mai\footnotemark[1]
	\and  
	Mikael Johansson{\thanks{Division of Decision and Control Systems, School of Electrical Engineering and Computer Science, KTH~Royal Institute of Technology, SE-100 44  Stockholm, Sweden. Emails: \tt\small\{maivv, mikaelj\}@kth.se.}}       
}

\begin{document}
\maketitle

\begin{abstract}

Stochastic gradient algorithms are often unstable when applied to functions that do not have Lipschitz-continuous and/or bounded gradients. Gradient clipping is a simple and effective technique to stabilize the training process for problems that are prone to the exploding gradient problem. Despite its widespread popularity, the convergence properties of the gradient clipping heuristic are poorly understood, especially for stochastic problems. This paper establishes both qualitative and quantitative convergence results of the clipped stochastic (sub)gradient method (SGD) for non-smooth convex functions with rapidly growing subgradients. Our analyses show that clipping enhances the stability of SGD and that the clipped SGD algorithm enjoys finite convergence rates in many cases. We also study the convergence of a clipped method with momentum, which includes clipped SGD as a special case, for weakly convex problems under standard assumptions. With a novel Lyapunov analysis, we show that the proposed method achieves the best-known rate for the considered class of problems, demonstrating the effectiveness of clipped methods also in this regime. Numerical results confirm our theoretical developments.

\end{abstract}

\section{Introduction}
We study stochastic optimization problems on the form
\begin{equation}
\label{eqn:objective}
		\underset{x\in \R^n}{\minimize} \,
%		\, f(x) 
%	\,\, \mbox{where} \,\, 
	 f(x) := \E_{P}[f(x;\statrv)]
	= \int_{\mc{S}} f(x; \statval) dP(\statval),
\end{equation}
where $\statrv\sim P$ is a random variable;  $f(x;\statval)$ is the instantaneous loss parameterized by $x$ on a sample $\statval\in\mc{S}$. Such problems are at the core of many  machine-learning applications, and are often solved using stochastic (sub)gradient methods.
%
%Stochastic (sub)gradient methods for solving problems  on the form~\eqref{eqn:objective} are horse in many  large-scale machine learning applications. 
In spite of their successes, stochastic gradient methods can be sensitive to their parameters \cite{NJLS09,AD19a} and have severe instability (unboundedness) problems when applied to functions that grow faster than quadratically in the decision vector $x$ \cite{And96, AD19a}. Consequently, a careful (and sometimes time-consuming) parameter tuning is often required for these methods to perform well in practice. Even so, a good parameter selection is  not sufficient to circumvent the instability issue on steep functions. 

Gradient clipping and the closely related gradient normalization technique are simple modifications to the underlying algorithm to control the step length that an update can make relative to the current iterate. These techniques enhance the stability of the optimization process, while adding essentially no extra cost to the original update. As a result,  gradient clipping has been a common choice in many applied domains of machine learning \cite{PMB13}. 

In this work, we consider gradient clipping applied to the classical SGD method. Throughout the paper, we frequently use the following clipping operator
\begin{align*}
	{\mathrm{clip}_\gamma: \R^n \to \R^n:  x \mapsto  \min\left\{1, \frac{\gamma}{\ltwo{x}}\right\}x,}
\end{align*}
which is nothing else but the orthogonal projection onto the $\gamma$-ball. 
It is important to note that for noiseless gradients, clipping just changes the magnitude and does not effect the search direction. However, in the stochastic setting, the expected value of the clipped stochastic gradient may point in a completely different direction than the true gradient. 

\noindent\textbf{Clipped SGD}\hspace*{0.5em} 
To solve problem \eqref{eqn:objective}, we use an iterative procedure that starts from $x_0\in\R^n$ and $g_0 \in \subdiff f(x_0, S_0)$ and generates a sequence of points $x_k\in\R^n$ by repeating the following steps for $k=0,1,2,\ldots$:
\begin{align}\label{alg:clip:sgd}
	x_{k+1} = x_k - \stepsize_k d_k, \quad 
	d_{k} = \tvclip{\gamma_k}{g_k}.
\end{align}
We refer to  $\stepsize_k$ as the $k$th stepsize and  $\gamma_k$ as the $k$th clipping threshold, while  $g_k=f'(x_{k}, \statrv_{k})$ is the $k$th stochastic subgradient or its mini-batch version, $g_k=\frac{1}{\batchsize_k}\sum_{i=1}^{\batchsize_k} f'(x_{k}, \statrv_{k}^i)$ if multiple samples are used in each iteration.

\subsection{Related work}
Our work is closely related to a number of topics which we briefly review below.

\paragraph{Gradient clipping} 
Gradient clipping and normalization were recognized early in the development of subgradient methods as a useful  tool to obtain convergence for rapidly growing convex functions \cite{Sho85, Erm88, AIS98}. For the normalized method, the seminal work~\cite{Sho85}  establishes a convergence rate for the quantity $\InP{g_k/\ltwo{g_k}}{x_k-x\opt}$  without any assumptions on $g_k$. By only requiring that subgradients are bounded on bounded sets, which always holds for continuous functions, the authors in \cite{AIS98} prove convergence in the objective value for the clipped subgradient method. 
{The work \cite{HLS15} considers normalized gradient methods for quasi-convex and locally smooth functions. }
Recently, the authors in \cite{ZHSJ19,ZJFW20} analyze clipped methods  for twice-differentiable functions satisfying a more relaxed condition than the traditional $L$-smoothness. However, much less is known in the stochastic setting. The work \cite{And96} proposes a method that uses two independent samples in each iteration and proves its almost sure convergence under the same growth condition used in \cite{AIS98}. The authors in {\cite{HLS15} establish certain complexity results for a mini-batch stochastic normalized gradient method under some strong assumptions on the closeness of all the generated iterates to the optimal solution and the boundedness of all the individual mini-batch functions. In \cite{ZHSJ19,ZJFW20},   stochastic clipped SGD methods are analyzed under the assumption that  the noise in the gradient estimate is bounded for every $x$ almost surely. Such assumption does not hold even for Gaussian noise in the gradients. } Finally, we refer to \cite{CM20, CSS19,GDG20} for recent theoretical developments for clipped and normalized methods on standard $L$-smooth problems.

\paragraph{Robustness and stability}  
The problems of robustness and stability in stochastic optimization  have been emphasized in many studies (see, e.g., \cite{NJLS09,  And96, AD19a, AD19b} and references therein). Much recent work on this topic concentrates around model-based algorithms that attempt to construct more accurate models of the objective than the linear one provided by the stochastic subgradient. 
When such models can be obtained and the resulting update steps can be performed efficiently, these methods often possess good stability properties and can be less sensitive to  parameter selection than traditional stochastic subgradient methods. For example, the work \cite{AD19a} establishes almost sure convergence of stochastic (approximate) proximal point methods under the arbitrary growth condition used in~\cite{AIS98, And96} and discussed above. In \cite{AD19b}, almost sure convergence of the so-called truncated method  is proven for convex functions that can grow polynomially from the the set of solutions.

\paragraph{Weakly convex minimization} 
The class of weakly convex functions is broad, allowing for both non-smooth and non-convex objectives, and has favorable structures for algorithmic foundations and complexity theory. Earlier works on weakly convex minimization \cite{Nur73,Rus87, EN98} establish qualitative convergence results for subgradient-based methods.
With the recent advances in statistical learning and signal processing, there has been an emerging line of work on this topic (see, e.g., \cite{DR18b,DG19,DD19}). Convergence properties have been analyzed for many popular stochastic algorithms such as: model-based methods \cite{DD19, DR18b, AD19b}; momentum extensions \cite{MJ20}; adaptive methods \cite{AMC20};  and more. 

\subsection{Contributions}
The performance of stochastic (sub)gradient methods depends heavily on how rapidly the underlying function is allowed to grow. Much convergence theory for these methods hinges on the $L$-smoothness assumption 
%\added{assumptions of} $L$-smooth gradients 
for differentiable functions or uniformly bounded subgradients for non-smooth ones. These conditions restrict the corresponding convergence rates to functions with at most quadratic and linear growth, respectively.
Beyond these well-behaved  classes, there is abundant evidence that SGD and its relatives may fail to converge. It is our goal in this work to show, both theoretically and empirically,  that the addition of the clipping step greatly improves the convergence properties of SGD. To that end, we make the following contributions:
\begin{itemize}[leftmargin=0.4cm]

\item We establish stability and convergence guarantees for clipped SGD on convex problems with arbitrary growth (exponential, super-exponential, etc.) of the subgradients. We show that clipping coupled with standard mini-batching suffices to guarantee almost sure convergence. Even more, a finite convergence rate can also be obtained in this setting.

\item  We then turn to convex functions with polynomial growth and show that without the need for mini-batching, clipped SGD can essentially achieve the same optimal convergence rate as for stochastic strongly convex and Lipschitz continuous functions.

\item  We consider a momentum extension of clipped SGD for weakly convex minimization under standard growth conditions. With a carefully constructed Lyapunov function, we are able to overcome the bias introduced by the clipping step and preserve the best-known sample complexity for this function class.

\end{itemize}

Our experiments on phase retrieval, absolute linear regression, and classification with neural networks reaffirm our theoretical findings that gradient clipping can: i) stabilize and guarantee convergence for problems with rapidly growing gradients; ii) retain and sometimes improve the best performance of their unclipped counterparts even on standard problems.  We note also that none of the convergence results in this work require hard-to-estimate parameters to set the clipping threshold. 

\paragraph{Notation} 
For any $x,y\in\R^n$, we denote by $\InP{x}{y}$ the Euclidean inner product of $x$ and $y$. 
We denote by $\subdiff f(x)$ the Fr\'{e}chet subdifferential of $f$ at $x$;  $f'(x)$ denotes any element of $\subdiff f(x)$. The $\ell_2$-norm is denoted by $\ltwo{\cdot}$. For a closed and convex set $\mc{X}$,  
the distance and the projection map are given respectively by:  
$\dist(x, \mc{X})=\min_{z\in\mc{X}}\ltwo{z-x}$ and ${\mathop{\Pi}}_{\mc{X}}(x)=\argmin_{z\in\mc{X}}\ltwo{z-x}$.  $\bindic{E}$ denotes the indicator function of an event $E$; i.e., $\bindic{E}=1$ if $E$ is true and 0 otherwise. The closed $\ell_2$-ball centered at $x$ with radius $r>0$ is denoted $\mathrm{B}(x, r)$.  We denote by $\mc{F}_k := \sigma(\statrv_0, \ldots, \statrv_{k-1})$ the $\sigma$-field formed by the first $k$ random variables $\statrv_0, \ldots, \statrv_{k-1}$, so that $x_k\in \mc{F}_k$. Finally,  we will impose the following basic assumption throughout the paper.
\begin{assumption}\label{assumption:unbiased}
Let $\statrv$ be a sample drawn from $P$ and $f'(x,\statrv) \in \subdiff f(x, \statrv)$, we have: $\E\left[f'(x, \statrv)\right] \in \subdiff f(x)$.
\end{assumption}

%\vspace{-0.2cm}
\section{Stability and its consequences for convex minimization}%\vspace{-0.1cm}
In this section, we study the stability of the clipped SGD algorithm and its consequence for the minimization of (possibly non-smooth) convex functions. 
We first specify the assumptions needed for the results in this section starting with the basic quadratic growth condition.
\begin{assumption}[Quadratic growth]\label{assumption:quadratic:growth}
 There exists a scalar $\mu>0$ such that
\begin{align*}
	 f(x) - f(x\opt) \geq \mu \dist\left(x, \Xopt\right)^2, \quad \forall x\in\dom(f).
\end{align*}
\end{assumption}%\vspace{-0.1cm}
Assumption~\ref{assumption:quadratic:growth} gives a \emph{lower bound} on the speed at which the objective $f$ grows away from the solution set $\Xopt$. Since we are interested in problems that may exhibit exploding subgradients, this growth condition is a rather natural assumption. Note also that in many machine learning applications, the addition of a quadratic regularization term to improve generalization results in problems which fundamentally have quadratic growth.
\begin{assumption}[Finite variance]\label{assumption:variance:cvx}
There exists a scalar $\sigma>0$ such that:
\begin{align*}
	\E\left[\ltwo{f'(x, S) - f'(x)}^2\right]  \leq  \sigma^2, \quad \forall x\in \dom(f),
\end{align*}
where  $f'(x) = \E\left[f'(x, S) \right]  \in \subdiff f(x)$.
\end{assumption} 
Finally, unless otherwise stated, we assume that the stepsizes $\stepsize_k$ are square summable but not summable:
\begin{align*}
	\stepsize_k \geq 0, 
	\quad
	\sum_{i=0}^{\infty} \stepsize_k = \infty,
	\quad \mbox{and} \quad
	\sum_{i=0}^{\infty} \stepsize_k^2 < \infty.
\end{align*}
Before detailing the stability and convergence analyses of  clipped SGD, Example~1  shows that even with stepsizes that are as small as $O(1/k)$, the vanilla SGD method may fail miserably when applied to a function satisfying Assumptions~\ref{assumption:quadratic:growth}--\ref{assumption:variance:cvx}. We refer to \cite{AD19a} for more examples of the potential instability of SGD.\medskip

\begin{example} [Super-Exponential Divergence of SGD]\label{example:sgd:divergence}
Let $f(x) = x^4/4 + \epsilon x^2/2$ with $\epsilon>0$ and consider the SGD algorithm applied to $f$ with the stepsizes  $\stepsize_k=\stepsize_1/k$:
\begin{align*}
	x_{k+1} = x_k - \frac{\stepsize_1}{k} \left(x_k^3 + \epsilon x_k \right).
\end{align*}
Then, if we let $x_1 \geq \sqrt{3/\stepsize_1}$, it holds for any $k\geq1$ that $\abs{x_k} \geq \abs{x_1} k!.$
\end{example}

Despite its simplicity, the example highlights that moving beyond standard (upper) quadratic models, SGD may fail to have any meaningful convergence guarantees.
 Our goal in this section is to: (i) show that with a simple clipping step added to SGD, the resulting algorithm becomes much more stable; and (ii) to prove  strong convergence guarantees for clipped SGD in new settings. Next, we state the first of these results:
\begin{proposition}[Stability]\label{propositiom:stability}
Let Assumptions~\ref{assumption:unbiased}, \ref{assumption:quadratic:growth}, and \ref{assumption:variance:cvx} hold. Let $\gamma_k \leq \gamma /\sqrt{\stepsize_k}$ for some $\gamma>0$. Let $C={\sigma^2}/(2\mu) + \gamma^2$, then, the iterates generated by the clipped SGD method satisfy
\begin{align}\label{eq:proposition:stability}
	\E\left[ \dist\left(x_{k}, \Xopt\right)^2 \right]  \leq \dist\left(x_{0}, \Xopt\right)^2 +
	 	C\sum_{i=0}^{k-1} \stepsize_i.
\end{align}
\end{proposition}%\vspace{-0.15cm}
Some remarks on Proposition~\ref{propositiom:stability} are in order. First, unlike SGD, where the distance to the optimal set may grow super-exponentially, the clipped version will not diverge faster than the sum of the used stepsizes. For example, with the stepsizes $O(1/k)$ in Example~\ref{example:sgd:divergence}, the sum is only of order $\log(k)$. {Such a guarantee will play a critical role in establishing all the convergence results in the subsequent sections.} Second, the proposition is reported for time-varying clipping thresholds to facilitate the proofs of some subsequent results. We note however that the similar estimate holds for the constant scheme  with a slightly different scaling constant. Finally, we mention that the bound \eqref{eq:proposition:stability} is similar to the classical results for the stochastic \emph{proximal point} iteration \cite[Theorem~6]{RB14}, but slightly weaker than the best bounds for that algorithm~\cite[Corollary~3.1]{AD19a}.

\subsection{Convergence under arbitrary growth}

Having studied the stability of clipped SGD, we now turn to its consequences for the actual convergence guarantees. We first remark that on \emph{deterministic} convex problems,  the procedure \eqref{alg:clip:sgd} is known to be convergent  under the very weak growth condition summarized in Assumption~\ref{assumption:Gbig} below~\cite{AIS98}. Concretely, the subgradients can grow arbitrarily (exponentially, super-exponentially, etc.) as long as they are bounded on bounded sets. However, the situation is less clear as stochastic noise enters the problem. Under \ref{assumption:Gbig},  similar convergence results have only been established for the stochastic proximal point method \cite{AD19a} {and a scaled stochastic approximation algorithm proposed in \cite{And96}.  Note that the former algorithm relies heavily on the ability to accurately model the objective and efficiently solve the resulting minimization problem in each iteration, while the later one needs two independent search directions to construct its upates.} Theorem~\ref{thrm:cvx:minibatch:asym} below demonstrates that gradient clipping coupled with mini-batching can also provide such a strong qualitative guarantee.

\begin{assumption}\label{assumption:Gbig}
There exits an increasing function $G_{\mathrm{big}}: \R_{+}\to [0,\infty)$ such that
\begin{align*}
	\E\left[\ltwo{f'(x, S)}^2\right]  
	 \leq G_{\mathrm{big}}(\dist\left(x, \Xopt\right)), \quad \forall x\in \dom(f).
\end{align*}
\end{assumption}  

\begin{theorem}\label{thrm:cvx:minibatch:asym}
Let Assumptions \ref{assumption:unbiased}, \ref{assumption:quadratic:growth}, and \ref{assumption:variance:cvx} hold. Let  $\gamma_k=\gamma$ for all $k$.  Consider for each $k$ a batch of samples $S_k^{1:m_k}$ and  let $x_k$ be generated  by the clipped SGD method with $g_k=\frac{1}{m_k}\sum_{i=1}^{m_k} f'(x, S_k^i)$.  Define $\varrho_k=\min\left\{1, \gamma/\ltwo{g_k}\right\}$ and $e_k = \dist\left(x_{k}, \Xopt\right)$, then
\begin{align*}
	\E\left[ e_{k+1}^2 \big | \mc{F}_{k} \right] 
%	\nonumber\\
%	 &\hspace{0.2cm}	
	 \leq
		\left(1 - \mu \stepsize_k 	\E\left[\varrho_k \big | \mc{F}_k \right] \right)
        e_k^2
		+ \frac{ \sigma^2\stepsize_k}{ \mu \batchsize_k} 
		+ \stepsize_k^2 \gamma^2.
\end{align*}
Suppose further that $\sum_{k=0}^\infty \stepsize_k/\batchsize_k < \infty$, then under Assumption~\ref{assumption:Gbig}, we have
%\begin{align*}
	$\dist\left(x_{k}, \Xopt\right) \lcas 0.$
%\end{align*}
\end{theorem}
Theorem~\ref{thrm:cvx:minibatch:asym} highlights the importance of the clipping step as no amount of samples in a batch can save SGD from divergence in this setting. In particular,  it implies that clipped SGD converges for any growth function provided that sufficiently accurate estimates of the subgradients can be obtained. 
This is in stark contrast to SGD without clipping, where, as Example~\ref{example:sgd:divergence} shows, the iterates may diverge  even in the noiseless setting when the objective function grows faster than the quadratic $x^2$. Since the stepsizes are square summable, taking $m_k=1/\stepsize_k$ suffices to guarantee $\sum_{k=0}^\infty \stepsize_k/\batchsize_k < \infty$.
%Again, we defer the detailed proof to Appendix~\ref{appendix:proof:thrm:cvx:minibatch:asym} and \replaced{only sketch the main arguments here. After having established the inequality in the theorem, we invoke}{sketch here only the main arguments. Suppose that we are given the inequality in the theorem, we can invoke} a well-known \replaced{result by Robbins and Siegmund on the almost sure convergence of nonnegative almost supermatingales}{Robbins-Siegmund almost supermartingale convergence lemma}\cite{RS71} to obtain 
%\begin{align*}
%	 \dist\left(x_{k}, \Xopt\right)^2 \lcas V_\infty < \infty 
%	\quad \mbox{and} \quad
%	\sum_{k=0}^{\infty} \stepsize_k \E\left[\varrho_k \big | \mc{F}_k \right]  \dist\left(x_{k}, \Xopt\right)^2 < \infty \quad \mbox{a.s.}
%\end{align*}
%We then use the first consequence to prove that the quantity $\E\left[\varrho_k \big | \mc{F}_k \right]$ is bounded away from zero almost surely. Therefore, the second consequence yields $\sum_{k=0}^{\infty} \stepsize_k \dist\left(x_{k}, \Xopt\right)^2 <\infty$ almost surely.  Since  $\sum_{k=0}^{\infty} \stepsize_k=\infty$, we deduce that $\liminf_{k\to \infty} \dist\left(x_{k}, \Xopt\right) = 0$ almost surely. But by the first consequence, the limit exists, and hence we must have $ \dist\left(x_{k}, \Xopt\right) \lcas 0$.

 It turns out that clipping can even provide finite convergence rate in this setting, as stated in the next result.
\begin{theorem}\label{thrm:cvx:minibatch:finite}
Let Assumptions \ref{assumption:unbiased}, \ref{assumption:quadratic:growth},  \ref{assumption:variance:cvx}, and \ref{assumption:Gbig} hold. Let $\stepsize_k=(k+1)^{-\tau}$ with $\tau\in (1/2,1)$ and let $x_k$ be generated by clipped SGD  using  batches of $m_k=1/\stepsize_k$ samples. 
Fix a failure probability $\delta\in (0,1)$, then for any $\epsilon>0$, there exists a numerical constant $c_0>0$ such that 
\begin{align*}
	\Pr\left(\dist\left(x_K, \Xopt\right)^2 \leq \epsilon \right) 
	\geq
		1 - \delta - \frac{\delta  \left({ \sigma^2  }/{\mu } +  \gamma^2\right) \sum_{k=0}^{K-1} \stepsize_k^2 }{e_0^2}  - \frac{c_0 \stepsize_K}{\epsilon}.
\end{align*} 
Furthermore, let $\varrho = \gamma/ (\gamma + G_{\mathrm{big}}^{1/2}(\dist\left(x_0, \Xopt\right)/\delta))$. If we take $\stepsize_k= \stepsize=\stepsize_0 K^{-\tau}$ with $\stepsize_0\leq1/(\mu\varrho)$,  and $\eta={\left( \sigma^2 /  \mu +  \gamma^2\right)}/{\mu \varrho}$,  then, for $K\in \N_+$ satisfying 
$ \mu\varrho\stepsize_0 K^{1-\tau} \geq \log\left({e_0^2 K^{\tau}}/{\eta \stepsize_0}\right)$, we have
\begin{align*}
	\dist\left(x_K, \Xopt\right)^2 \leq \frac{2 \eta \stepsize_0}{\delta K^{\tau}},
\end{align*}
with probability at least
$1 - 2\delta -  \delta\cdot\frac{ \left({ \sigma^2  }/{\mu } +  \gamma^2\right) \stepsize_0^2}{\dist\left(x_0, \Xopt\right)^2 K^{2\tau-1}}$.
\end{theorem}
The first result in the theorem refines the asymptotic guarantee in Theorem~\ref{thrm:cvx:minibatch:asym} for general time-varying stepsizes and the second one shows the iteration complexity for a constant stepsize and fixed mini-batch size. We have the following remarks: (i) Setting $\tau$ close to one in the second claim yields a bound with a similar order-dependence on K and $\delta$ as for strongly convex and Lipschitz continuous $f$ \cite[eq.~(4.2.61)]{Lan20}. Note that the last term in the last probability bound is negligible for large $K$;  (ii) The proof of the theorem is motivated by a technique developed in \cite[Lemma~3.3]{DDC19} to bound the escape probability of their algorithm's iterates.

\subsection{Convergence under polynomial growth}

For the final set of theoretical results of the section, we consider a more specific function class for which we derive the convergence rate of clipped SGD without the need \replaced{for}{of} mini-batching. In particular, we impose the following conditions on the stochastic subgradients.
\begin{assumption}\label{assumption:polynomial:growth}
 There exist real numbers $L_0, L_1, \sigma \geq 0$ and $ 2 \leq p < \infty$ such that for all $x\in\dom(f)$:
\begin{align*}
	&\E\left[\ltwo{f'(x, S)}^2\right]  \leq  L_0 + L_1\dist\left(x, \Xopt\right)^{2(p-1)}, 
	\\
	% \label{assumption:central:moment}
	&\E\left[\ltwo{f'(x, S) - f'(x)}^{2(p-1)}\right]  \leq  \sigma^p, 
\end{align*}
 where  $f'(x) = \E\left[f'(x, S) \right]  \in \subdiff f(x)$.
\end{assumption}
Note that when $p=2$, we have the standard (upper) quadratic growth model \cite{PJ92}. For general values of $p$,  Assumption~\ref{assumption:polynomial:growth} implies that 
\begin{align*}
	\ltwo{f'(x)}^2
%	=  \ltwo{\E\left[f'(x, S)\right]}^2
	\leq \E\left[\ltwo{f'(x, S)}^2 \right]
	\leq  L_0 + L_1\dist\left(x, \Xopt\right)^{2(p-1)}
\end{align*}
which, since $f$ is assumed to be convex, guarantees that 
%We can thus deduce from  the convexity of $f$ that
\begin{align*}
	f(x) - f(x\opt) 
%	\leq \InP{f'(x)}{x-x\opt} 
	\leq \sqrt{L_0} \dist\left(x, \Xopt\right) + \sqrt{L_1}\dist\left(x, \Xopt\right)^{p}.
\end{align*}
We thus allow the function $f$ to grow polynomially from the set of optimal solutions.
For example, $f(x)=x^4/4 + \epsilon x^2/2$ satisfies the assumption with $L_0=L_1=2(1+\epsilon)$ and $p=4$. The second condition in~\ref{assumption:polynomial:growth} requires that the $2(p-1)$th central moment is bounded, which amounts to finite variance when $p=2$. 
We mention that a closely related assumption has been used in \cite[Assumption~A3]{AD19b} to analyze a method analogous to the classical Polyak subgradient algorithm. The only difference is  in the second condition, where they require bounded variation of a quantity involving the objectives instead of the subgradients. This is because $f(x,S)$ and $f(x\opt)$ are used to construct their updates.

The next lemma explicitly bounds the expected norm of the subgradients and the distance between the iterates and the set of optimal solutions.  The proof of this lemma follows the same arguments in \cite[Lemma~B2]{AD19b} and is reported in 
%Appendix~E
Appendix~\ref{appendix:proof:lem:xdiff:grad:bound} 
for completeness.
\begin{lemma}\label{lem:xdiff:grad:bound}
{Let Assumptions \ref{assumption:unbiased}, \ref{assumption:quadratic:growth}, and \ref{assumption:polynomial:growth} hold. Let $x_k$ be generated by clipped SGD using $\stepsize_k = \stepsize_0 (k+1)^{-\tau}$ with $ \tau \in (1/2,1)$, then there exist positive real constants $D_0, D_1, G_0, G_1$ (independent of $k$) such that
\begin{align*}
		\E\left[\ltwo{f'(x_k, S)}^2\right]   &\leq G_0 + G_1 k^{(p-1)(1-\tau)},
		\\
		\E\left[\dist\left(x_k, \Xopt\right)^{4(p-1)} \right] &\leq D_0 + D_1 k^{2(p-1)(1-\tau)}.
\end{align*}
More specifically, if we define for $q\geq 2$ the quantities:
\begin{align*}
	P_0(q) &:= 2^{\frac{q}{2}} \dist(x_0, \Xopt)^q\quad \mbox{and} \quad
	\nonumber\\
	P_1(q) &:= 	\left(  (2\gamma)^q + \mu^{-\frac{q}{2}} 	\sigma^{\frac{q}{4}+1} \right) \left(\frac{2\stepsize_0}{1-\tau}\right)^{\frac{q}{2}},
\end{align*}
then $G_0, G_1, D_0, D_1$ are given explicitly by
\begin{align*}
	\begin{array}{ll}
	 G_0 = L_0 + L_1 P_0(2(p-1)),  &
	 G_1 = L_1 P_1(2(p-1)), 
	 \\
	 D_0 = P_0(4(p-1)),  &  
	 D_1 = P_1(4(p-1)).
	\end{array} 	
	\end{align*}}
\end{lemma}
The lemma reveals an attractive property: the subgradients at the iterates can be made small by setting $\tau$ close to one, no matter the value of $p$. This brings us to a position close to where we would have been if we had assumed  Lipschitz continuity of $f$ in the first place. {The difference is, however, that the preceding guarantees hold w.r.t the full expectation while in the alternative case, one is given a priori an upper-bound on the quantity $\E[\ltwo{f'(x_k, S)}^2 \big | \mc{F}_k]$.} We can now state the main result of this subsection. 
\begin{theorem}\label{thrm:cvx:finite:tvclip}
Let Assumptions \ref{assumption:unbiased}, \ref{assumption:quadratic:growth},  and \ref{assumption:polynomial:growth} hold. Let $x_k$ be generated by clipped SGD using $\stepsize_k = \stepsize_0 (k+1)^{-\tau}$ with $ \tau \in (1/2,1)$ and $\gamma_k=\gamma/\sqrt{\stepsize_k}$, then there exists a numerical constant $C$ such that
\begin{align*}
	\E \left[\dist\left(x_{k+1}, \Xopt\right) ^2 \right] 
	\leq
			\left(1- \frac{\mu\stepsize_0}{(k+1)^\tau} \right) \E \left[\dist\left(x_k, \Xopt\right)^2 \right]
%	\nonumber\\
%	&\hspace{0.45cm}
	+ \frac{C}{(k+1)^{2 \left(1- p(1-\tau)\right)}}.
\end{align*}
Furthermore, we we take $\tau = 1- \epsilon$ for some $\epsilon>0$, then 
\begin{align*}
	\E \left[\dist\left(x_k, \Xopt\right) ^2 \right] 
	\leq
			\frac{C}{\mu\stepsize_0}\frac{1}{k^{1+\epsilon(1 - 2 p)}}
			+ o\left(\frac{1}{k^{1+\epsilon (1- 2p})}\right).
\end{align*}

%\begin{align*}
%	\E \left[\ltwo{x_{k+1}-x\opt} ^2 \right] 
%	&\leq
%			\left(1- \frac{\mu\stepsize_0}{(k+1)^\tau}\right) \E \left[\ltwo{x_k-x\opt} ^2 \right]
%			+ \frac{C}{(k+1)^{2 \left(1- p(1-\tau)\right)}}.
%\end{align*}
%Furthermore, we we take $\tau = 1- \epsilon$ for some $\epsilon>0$, then 
%\begin{align*}
%	\E \left[\ltwo{x_{k+1}-x\opt} ^2 \right] 
%	&\leq
%			\frac{C}{\mu\stepsize_0}\frac{1}{k^{1+\epsilon(1 - 2 p)}}
%			+ o\left(\frac{1}{k^{1+\epsilon (1- 2p})}\right).
%\end{align*}
\end{theorem}
Some remarks on Theorem~\ref{thrm:cvx:finite:tvclip} are in order: 
\begin{enumerate}[label=(\roman*)]
\item {The numerical constant $C$ can be computed  as
\begin{align*}
	C = (2\gamma^2/\mu)  (L_0^2 + L_1^2  (D_0 + D_1))
 + G_0 + G_1, 
\end{align*}
where $D_0$, $D_1$, $G_0$, and $G_1$ are given in Lemma~\ref{lem:xdiff:grad:bound}. If $f$ is Lipschitz continuous ($L_1=0$), then $C$ reduces to $C=  (2\gamma^2/\mu)  L_0^2 + L_0$. Thus, by setting $\gamma = O(\sqrt{\mu/L_0})$ so that $C=O(L_0)$, we recover the similar order-dependence on $L_0$ as in the standard bound for unclipped SGD \cite{Lan20}.}

\item{ Despite the polynomial growth condition,  the first bound in the theorem is quite close  to the classical estimate for SGD ($\tau=1$), when applied to Lipschitz continuous $f$ using the \emph{small} stepsize $O(1/k)$ \cite{Lan20, Bec17, NJLS09}. Note also that our guarantee is valid for a wide range of \emph{large} stepsizes $\stepsize_k = O(k^{-\tau})$ with $\tau \in (1/2,1)$, which are more robust than the classical $O(1/k)$ \cite{NJLS09}.}

\item { The convergence result follows from a direct application of Chung's lemma~\cite[Lemma~4]{Chu54} to the first inequality in Theorem~\ref{thrm:cvx:finite:tvclip}. The little $o$ term in the estimate vanishes exponentially fast with the sum of the used stepsizes in the manner of \cite[p. 467]{Chu54} (see also \cite{BM11}). 
Since we are free to pick $\epsilon>0$, we  can, in principle, guarantee
a rate that is arbitrarily close to $O(1/ k)$. Recall that $O(1/ k)$ is also the  optimal convergence rate for (stochastic) strongly convex and Lipschitz continuous functions (two contradicting conditions) \cite[Section~2.1]{NJLS09}. 
\replaced{Hence, gradient clipping is able to essentially preserve}{Also here, we see that gradient clipping can almost preserve} the optimal rate of SGD while \replaced{also supporting}{allowing to support} a broad class of functions on which SGD and its relatives would diverge super-exponentially.}

\end{enumerate}

\section{Non-asymptotic convergence for weakly convex functions}
{The previous section demonstrates that gradient clipping can greatly improve the performance of SGD when applied to convex functions with rapidly growing subgradients. We now turn to non-asymptotic convergence analysis of clipped methods under standard growth conditions, but for a much wider class of weakly convex functions. Our goal is to show that the sample complexity of the clipped methods matches the best-known result for weakly convex problems, emphasizing the effectiveness of gradient clipping for a wide range of problem classes.}

Recall that that $f: \R^n \to \R\cup\{+\infty\}$ is called $\rho$-weakly convex if  $f + \frac{\rho}{2}\ltwo{\cdot}^2$ is convex. Such functions satisfy the following inequality for any $x, y\in \R^n$ with $g\in \subdiff f(x)$:
\begin{align*}
	f(y) \geq f(x) + \InP{g}{y-x} - (\rho/2) \norm{y-x}_2^2. 
\end{align*}
\replaced{Weakly convex optimization problems arise naturally in applications described by compositions}{
A major source of weak convexity comes from the composition} of the form
$f(x) = h(c(x)),$ where $h: \R^m \to \R$ is convex and $L_h$-Lipschitz and $c: \R^n \to \R^m$ is a smooth map with $L_c$-Lipschitz Jacobian. Note also that all convex functions and \added{all} differentiable functions with Lipschitz continuous gradient are weakly convex. We refer to \cite{DR18b, DD19, AD19a, MJ20} for practical applications as well as recent theoretical and algorithmic developments \replaced{for}{of} this function class.

\medskip

\noindent\textbf{Algorithm}\vspace*{0.5em} We consider the following momentum extension of clipped SGD:
\begin{subequations}
\label{alg:clip:shb}
\begin{align}
	x_{k+1} &= x_k - \stepsize_k d_k\\
	d_{k+1} &= \clip{(1-\mmt_k) d_k + \mmt_k g_{k+1}}.
\end{align}
\end{subequations}
Here, $g_{k+1}=f'(x_{k+1}, \statrv_{k+1})$, $\stepsize_k>0$ is the stepsize, $\mmt_k\in (0,1]$ is the momentum parameter, and $\thres>0$ is the clipping threshold. The algorithm is initialized from $x_0\in \R^n$ and $d_0=\clip{g_0}$ with $g_0\in \subdiff f(x_0, S_0)$, and generates the sequences $x_k \in \R^n$ (iterates) and $d_k\in R^n$ (search directions).
This algorithm goes back to at least \cite{GB72}, and in the sequel, we term procedure~\eqref{alg:clip:shb} as clipped stochastic heavy ball (SHB). 

Next we state the standing assumption in this section. 
\begin{assumption}\label{assumption:wcvx:gradient:boundedness} 
There exists a  positive real constant $L$ such that
\begin{align*}
\E \left[\ltwo{ f'(x, S)}^2\right] \leq L^2, \quad \forall x \in \dom(f).
\end{align*}
\end{assumption}
This is a very basic assumption in non-smooth optimization \cite{NJLS09,DD19}.\medskip

As the function $f$ is neither smooth nor convex, even measuring the progress to a stationary point for $f$ is a challenging task. A common practice is then to use the norm of the gradient of the Moreau envelope as a proxy for near-stationarity \cite{DD19}.  This is possible since weakly convex functions  admit an implicit smooth approximation through the classical Moreau envelope:
\begin{align}\label{eq:moreau:env}
	f_\lambda(x) = \inf_{y\in\R^n} \left\{f(y)+ \frac{1}{2\lambda} \ltwo{x-y}^2\right\}.
\end{align}
For $\lambda<\rho^{-1}$, the point achieving $f_\lambda(x)$ in \eqref{eq:moreau:env}, denoted by $\prox{\lambda f}{x} $,  is unique and given by:
\begin{align}
	\label{eq:prox}
	\prox{\lambda f}{x} = \argmin_{y\in\R^n} \left\{f(y) +  \frac{1}{2\lambda}\ltwo{x-y}^2\right\}.
\end{align}
With these definitions, for any $x\in\R^n$, the point $\hat{x}=\prox{\lambda f}{x}$ satisfies: 
\begin{align}\label{eq:stationary:measure}
\begin{cases}
	\ltwo{x-\hat{x}} =  \lambda \ltwo{\grad{f_{\lambda}(x)}},\\
	\dist(0, \subdiff f(\hat{x})) \leq \ltwo{\grad{f_{\lambda}(x)}}.
\end{cases}
\end{align}
Thus, a small gradient $\ltwo{\grad{f_{\lambda}(x)}}$ implies that $x$ is close to a point $\hat{x}$ that is near-stationary for $f$. 

The lemma below summarizes two useful properties of the Moreau envelope that will be used frequently in our convergence analysis \cite{HL93, DD19}.
\begin{lemma}[Moreau envelope]\label{lem:moreau:env}
Suppose that $f: \R^n \to \R\cup\{+\infty\}$ is a $\rho$-weakly convex function. Then, for a fixed parameter $\lambda^{-1} > 2\rho$, the following hold:
\begin{description}[leftmargin=0.45cm]
\item[1.]  $f_\lambda$ is continuously differentiable with the gradient given by
\begin{align*}
	\grad{f_{\lambda}(x)} = \lambda^{-1}\left(x-\prox{\lambda f}{x} \right).
\end{align*}
\item[2.]  $f_\lambda$ is $(1/\lambda)$--smooth, i.e., for all $x, y\in \R^n$:
\begin{align*}
			\ltwo{\grad{f_{\lambda}(x)} - \grad{f_{\lambda}(y)}} &\leq \frac{1}{\lambda}\ltwo{x-y}
			\nonumber\\
			 \big| f_\lambda(y) - f_\lambda(x) - \InP{\grad{f_{\lambda}(x)}}{y-x} \big|
		&\leq 		\frac{1}{2\lambda}\ltwo{x-y}^2.
\end{align*}
\end{description}
\end{lemma}
As for most convergence analyses of subgradient-based methods, we aim to establish  the following \added{per-iterate} estimate (see, e.g., \cite{NJLS09,DD19,GL13}):
\begin{align}\label{eq:stoc:analysis}
	\E[V_{k+1}] \leq \E[V_k] - c_0\alpha_k\, \E[e_k] + c_1\alpha_k^2.
\end{align}
Here $e_k$ is some stationarity measure, $V_k$ is a Lyapunov function, $\alpha_k$ is the stepsize, and $c_0, c_1$ are some real constants. As discussed above, for minimization of weakly convex functions  it is natural to consider $e_k=\Vert \nabla f_{\lambda}(x_k)\Vert^2$. It now remains to find an appropriate Lyapunov function $V_k$.
To build up our $V_k$, we will go through a number of supporting lemmas. We begin with the one that  concerns the search direction $d_k$.
\begin{lemma}\label{lem:xdiff:non-smooth}
Let Assumptions \ref{assumption:unbiased} and \ref{assumption:wcvx:gradient:boundedness} hold. Let $\mmt_k= \nu\stepsize_k$ for some constant $\nu>0$ such that $\mmt_k\in (0,1]$. Let $x_k$ be generated by  the clipped SHB method, then
\begin{align}\label{eq:lem:xdiff:non-smooth}
	\hspace{-0.2cm}	f(x_{k}) 
	+	
	\E\left[		
		\frac{1-\mmt_k}{2\nu}
		\ltwo{d_{k}}^2 \bigg|\mc{F}_{k}	\right]		
	&\leq
	f(x_{k-1})
	+
	\frac{1-\mmt_{k-1}}{2\nu}
	\ltwo{d_{k-1}}^2
	\nonumber\\
	&\hspace{0.45cm}
		- \stepsize_k\,\E\left[\ltwo{d_k}^2 \big|\mc{F}_{k}	\right]	
		+
		\frac{\stepsize_{k-1}^2}{2}
		\left(\frac{\nu L^2}{1-\mmt_{0}}+\rho\gamma^2\right).
\end{align}
\end{lemma}
It is interesting to note that, despite the bias introduced by the clipping operator, the estimate in \eqref{eq:lem:xdiff:non-smooth} is equivalent to that of in \cite[Lemma~3.1]{MJ20}.  Moreover, our new proof is \added{arguably} simpler, more intuitive and applicable to both constant and time-varying parameters.

The next lemma brings the gradient of the Moreau envelope to the stage.
\begin{lemma}\label{lem:wcvx:Wfunc}
Let Assumptions \ref{assumption:unbiased} and \ref{assumption:wcvx:gradient:boundedness} hold. Let $\mmt_k= \nu\stepsize_k$ for some constant $\nu>0$ such that $\mmt_k\in (0,1]$. Let $x_k$ be generated by  clipped SHB with $\gamma \geq 2L$ and define 
\begin{align*}
	W_{k} = \frac{1}{2\nu}\ltwo{d_{k} - \grad{f_{\lambda}(x_{k})}}^2 - \frac{1}{2\nu}\ltwo{\grad{f_{\lambda}(x_{k})}}^2
			+ f(x_{k}).
\end{align*}
Let $C= \nu L^2 + \rho\thres^2/2 $, then for any $k\in\N$, we have
\begin{align}\label{eq:lem:Wk}
	\E\left[W_{k}\big | \mc{F}_k\right]
	&\leq
		W_{k-1}
		-\stepsize_{k-1} 	\E\left[ \InP{g_k}{\grad{f_{\lambda}(x_k)}} \big | \mc{F}_k\right] 
	\nonumber\\
	&\hspace{0.45cm}
		+\stepsize_{k-1} \InP{d_{k-1}}{\grad{f_{\lambda}(x_{k-1})}}
		+ \frac{\stepsize_{k-1}}{\lambda\nu } \ltwo{d_{k-1}}^2
		+ C\stepsize_{k-1}^2.
\end{align}
\end{lemma}%\vspace{-0.2cm}
To motivate the introduction of $W_k$, we take a step back and consider the problem where $f$ is assumed to be $L$-smooth and \replaced{no gradient clipping is applied}{with no gradient clipping}. In this case, procedure~\eqref{alg:clip:shb} is identical to the algorithm \added{which was} analyzed in  \cite{RS83}\replaced{ using a Lyapunov function on the form}{. The authors introduced there the Lyapunov function of the form}
\begin{align*}
	V_k = \nu f(x_k) + \half\ltwo{d_k - \grad{f(x_k)}}^2 + \half\ltwo{d_k}^2.
\end{align*}
The key insight here is to view the sequence of directions $d_k$ as estimates of the true gradients $\grad{f(x_k)}$. With reasonable assumptions, the term $\E[\ltwo{d_k - \grad{f(x_k)}}^2]$ can indeed be driven to zero \cite[Theorem~1]{RS83}. Since our $f$ is non-smooth, this approach is not immediately applicable to our problem. Nevertheless,  we observe that it is useful to view $d_k$ as an estimate of the gradient of the Moreau envelope. This is the reason why the term $\ltwo{d_{k} - \grad{f_{\lambda}(x_{k})}}^2$ appears in  $W_k$, while the other terms arise from the algebraic manipulations to satisfy \eqref{eq:lem:Wk}. Finally, due to the presence of the clipping step, some extra care is needed to make the intuition work. In particular, since the $d_k$'s always belong to $\mathrm{B}(0,\gamma)$, we cannot expect that they approximate the $\grad{f_{\lambda}(x_k)}$ unless these also belong to the $\gamma$-ball. It turns out that  setting $\gamma\geq 2L$ suffices to ensure that $\grad{f_{\lambda}(x_k)} \in \mathrm{B}(0,\gamma)$.  

We now have all the ingredients needed to construct the ultimate Lyapunov function:
\begin{lemma}\label{lem:Vfunc:non-smooth}
Assume the same setting of Lemma~\ref{lem:wcvx:Wfunc}. Let $\lambda>0$ be such that $\lambda^{-1}\geq 2\rho$ and consider the function:
\begin{align*}
	V_{k} 
	= 
		f_\lambda(x_k) + W_k
		+
		\frac{f(x_k)}{\lambda\nu} 
		+ \left(\frac{1-\mmt_{k}}{2\lambda\nu^2} + \frac{\stepsize_k}{\lambda\nu}\right)\ltwo{d_k}^2.
\end{align*}
Then, for any $k\in \N_{+}$,
\begin{align}\label{eq:lem:Vfunc:non-smooth}
	 \E\left[V_{k}|\mc{F}_{k}\right]		 
	\leq
		 V_{k-1}	
		-
		\frac{\stepsize_{k-1}}{2}		
		\ltwo{\grad{f_{\lambda}(x_k)}}^2
		+
		C \stepsize_{k-1}^2,
\end{align}
where $C =  \lambda^{-1}\gamma^2 (1+ \rho/(2\nu)) + \nu L^2(1+ 1/(2\lambda\nu (1-\mmt_0)))$. 
\end{lemma}

Finally, the following complexity result \replaced{follows by a standard argument}{ is a standard deduction} from   \eqref{eq:lem:Vfunc:non-smooth}.
\begin{theorem}\label{thrm:wcvx}
Let Assumptions \ref{assumption:unbiased} and \ref{assumption:wcvx:gradient:boundedness} hold. Let $k^*$ be sampled randomly from $\{0,\ldots,K-1\}$ with $\Pr(k^*=k+1)= \stepsize_k/\sum_{i=0}^{K-1}\stepsize_i$. Let $\Delta=f(x_0)-\inf_{x}f(x)$ and let $C$ be given in \eqref{eq:lem:Vfunc:non-smooth}. Then, under the same setting of Lemma~\ref{lem:Vfunc:non-smooth}, we have
\begin{align*}
	\E\left[\ltwo{\grad{F_{\lambda}(x_{k^*})}}^2\right]	
	\leq
		2\cdot
		\frac{		
			\xi \Delta
			+ 2L^2/\nu 
			+	C \sum_{i=0}^{K-1} \stepsize_i^2		
		}{
			 \sum_{i=0}^{K-1} \stepsize_i
		},	
\end{align*}
where $\xi= 2+1/(\lambda\nu)$.
Furthermore, if we set $\stepsize={\stepsize_0}/{\sqrt{K}}$ and $\nu=1/\stepsize_0$ for some real $\stepsize_0>0$
\begin{align*}
	\E\left[\ltwo{\grad{F_{1/(2\rho)}(x_{k^*})}}^2\right]
	\leq 
		2\cdot
	\frac{
	   \xi \Delta + 2L^2/\nu+ C\stepsize_0^2
	}{
		\stepsize_0\sqrt{K}
	}.	
\end{align*}
Finally,  if $\stepsize_0$ is set to $1/\rho$ and $K\geq 2$, we obtain 
\begin{align*}
	\E\left[\ltwo{\grad{f_{1/(2\rho)}(x_{k^*})}}^2\right]
	\leq 
		8 \cdot
	\frac{
		\rho\Delta
		+ 
		\gamma^2
	}{
		\sqrt{K}
}.
\end{align*}
\end{theorem}
The rate achieved by the clipped SHB is of the same order as the best-known result for stochastic weakly convex   problems (see, e.g., \cite[Theorem~1]{DD19} and \cite[Theorem~1]{MJ20}). By inspection, all the proofs and convergence results in this section can also be extended (often with significant simplifications)  to the case of clipped SGD. The choice $\nu=1/\stepsize_0$ is just for simplicity; we \replaced{can choose any value of}{are free to pick} $\nu$ as long as $\mmt_k=\nu\stepsize_k\in(0,1]$. Since $\mmt_k=\nu \stepsize_k=O(1/\sqrt{k})$, one can put much more weight on the momentum term than on the fresh subgradient in the search directions $d_k$. 
As both $\stepsize_k$ and $\mmt_k$ have the same scale, the algorithm can be seen as a single time-scale method \cite{GRW18, RS83}. 

\section{Experimental results}

For the first two problems, we set up our experiments as follows. We fix  $m=500$,  $n=50$ and generate $A\in \R^{m\times n}$ as $A=Q D$, where $Q$ is a matrix with standard normal distributed entries, and $D$ is a diagonal matrix with linearly spaced elements between $1/\kappa$ and $1$. Here,  $\kappa\geq 1$ represents a condition number \replaced{which we}{and is} set to $\kappa=10$ in all experiments. The algorithms are all randomly initialized at $x_0\sim \mc{N}(0, 1)$ \replaced{and we use the stepsize}{. We set the stepsize as} $\stepsize_k=\stepsize_0 (k+1)^{-1/2}$, where $\stepsize_0$ is an initial stepsize. We also refer \added{to} $m$ stochastic iterations as one epoch (pass over \added{the} data). Within each individual run, we set the maximum number of epochs to 500. \replaced{Each plot reports the results of 30 experiments, visualized as the}{We perform for each plot 30 experiments and report the} median of the quantity of interest and the corresponding $90\%$ confidence interval.  Finally, the so-called epoch-to-$\epsilon$-accuracy is defined as the smallest number of epochs $q$ needed to reach $f(x_{m\cdot  q}) -  f(x^\star) \leq \epsilon$.

\begin{figure*}[t!]
	\begin{minipage}{0.495\textwidth}
	\centering
	{\includegraphics[width=1.\textwidth]{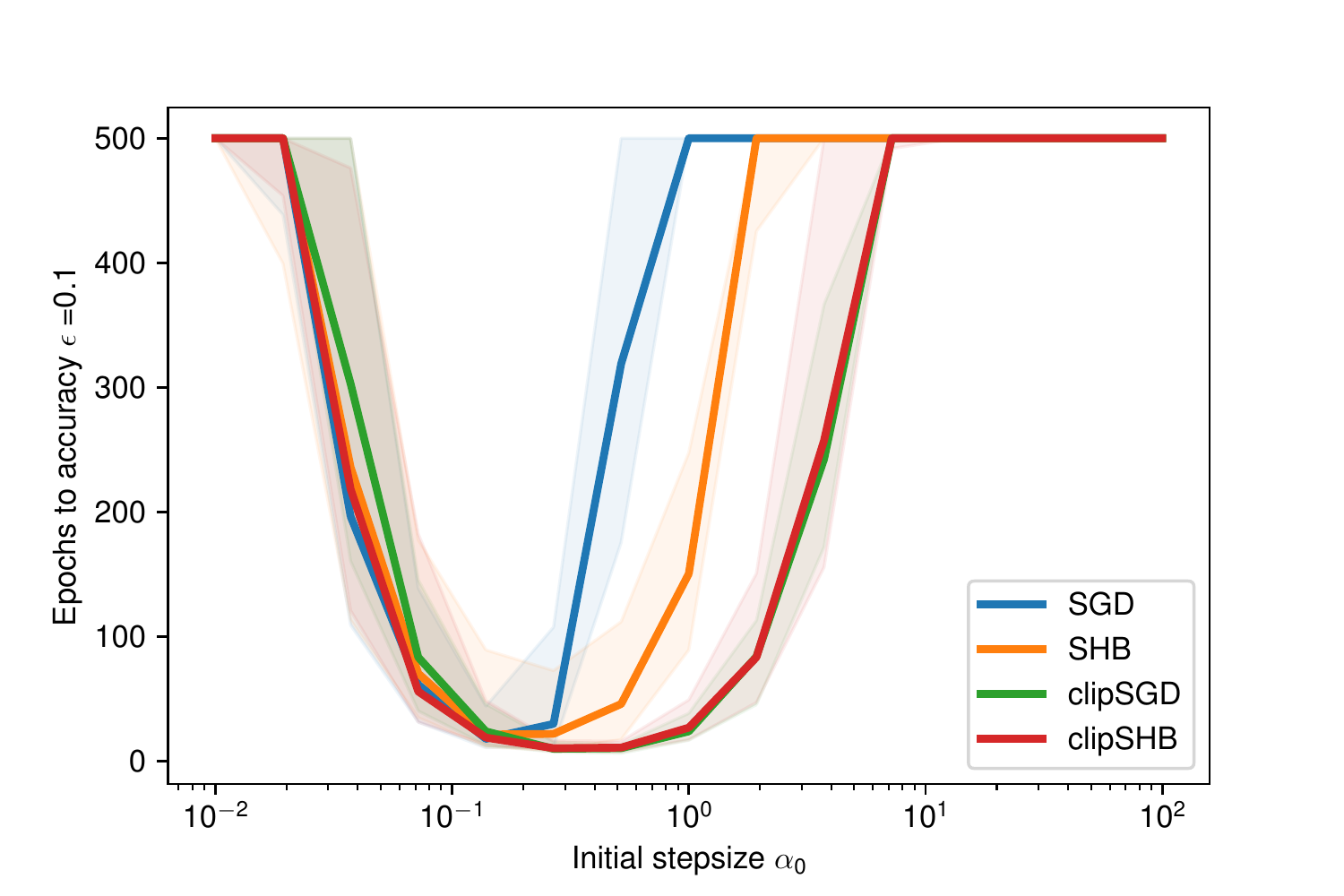}}
	\subcaption{ $1-\mmt = 0.9$, $\epsilon=0.1$}
	\end{minipage}
%	\hskip -0.15in
	\begin{minipage}{0.495\textwidth}
		\centering
		{\includegraphics[width=1.\textwidth]{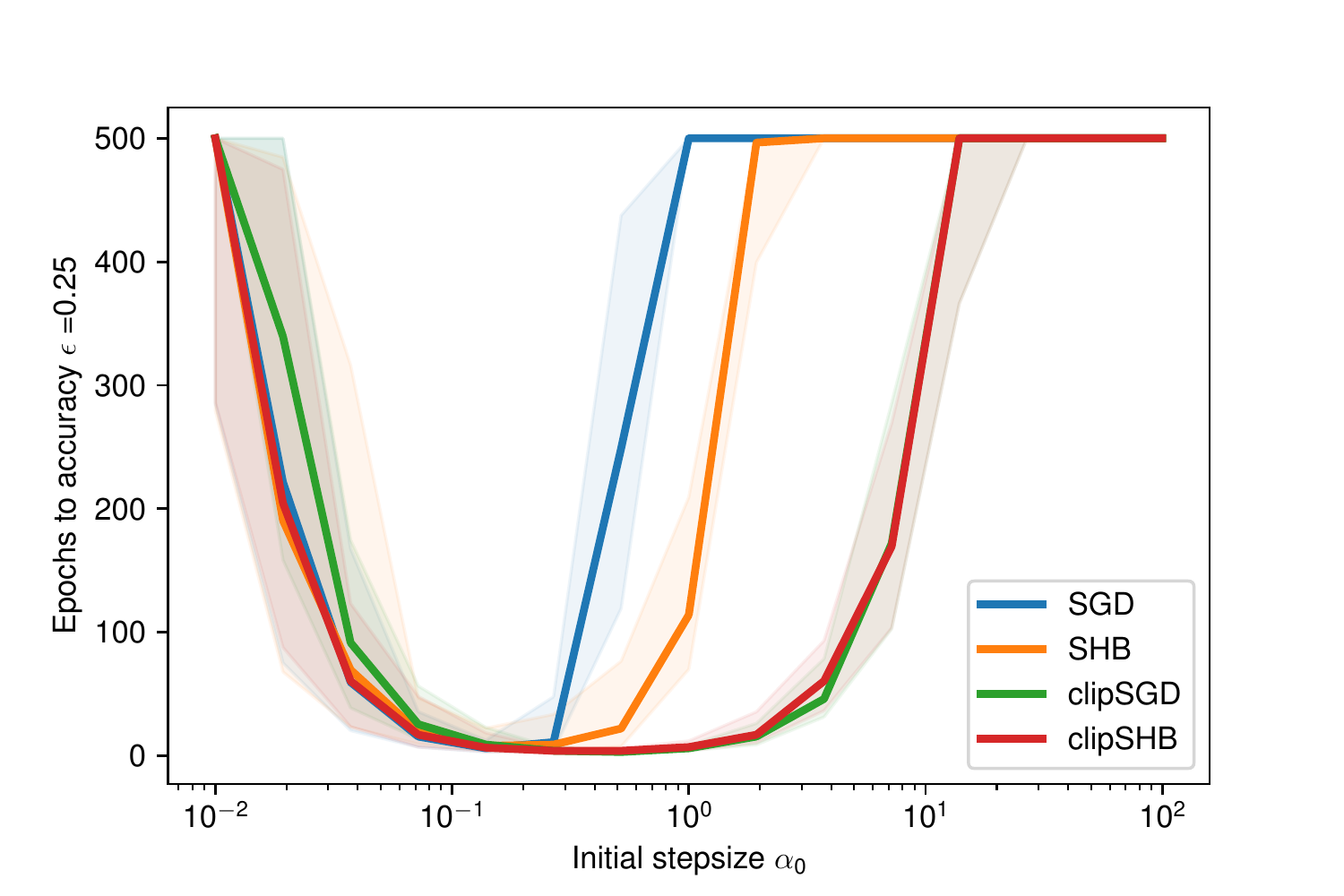}}
		\subcaption{ $1-\mmt = 0.9$, $\epsilon=0.25$}	
	\end{minipage}	
%	\vskip -0.1cm
	\centering
	\begin{minipage}{0.495\textwidth}
		\centering		{\includegraphics[width=1.\textwidth]{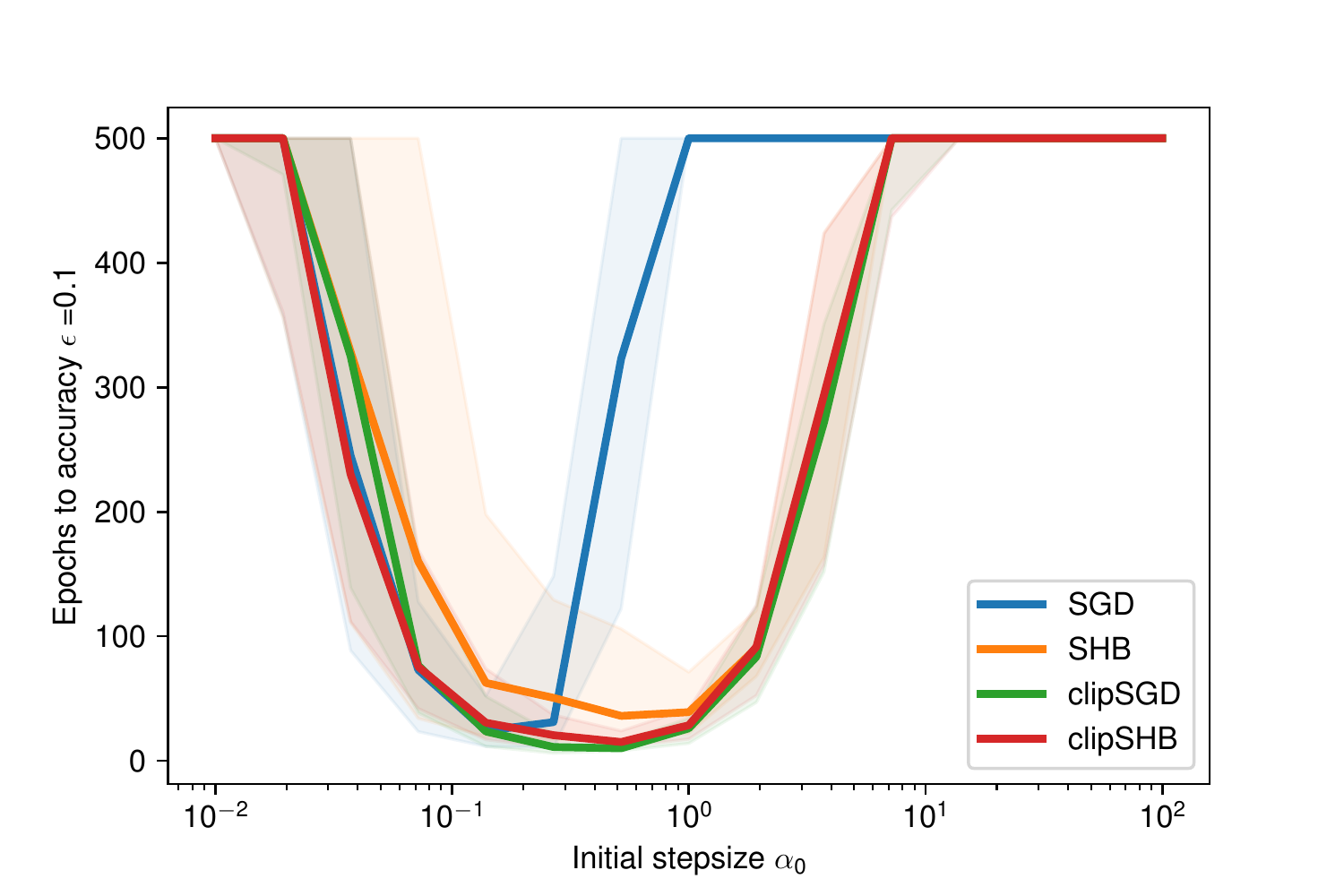}}
		\subcaption{$1-\mmt = 0.99$, $\epsilon=0.1$}	
		\end{minipage}
%		\hskip -0.15in
	\begin{minipage}{0.495\textwidth}
		\centering 
		{\includegraphics[width=1.\textwidth]{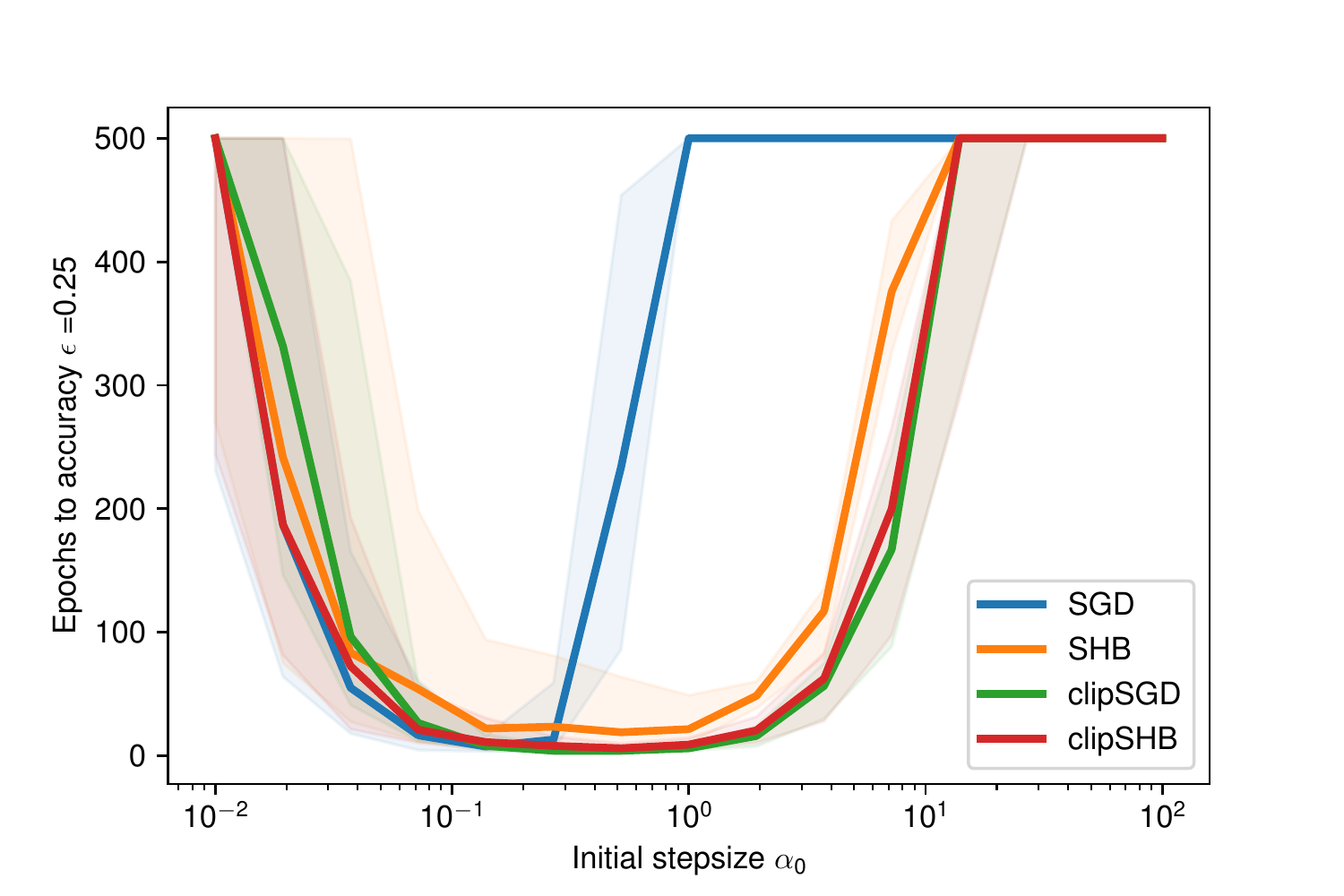}}
		\subcaption{ $1-\mmt = 0.99$, $\epsilon=0.25$}	
	\end{minipage}	
%	\vskip -0.1cm
	\caption{The number of epochs to achieve $\epsilon$-accuracy versus initial stepsize $\stepsize_0$ for phase retrieval with $\gamma=10$.}\label{fig:PR:epoch2eps}
\end{figure*}
\begin{figure*}[h!]
	%\vskip -0.15in
	\centering
	\begin{minipage}{0.495\textwidth}
		\centering		{\includegraphics[width=1.\textwidth]{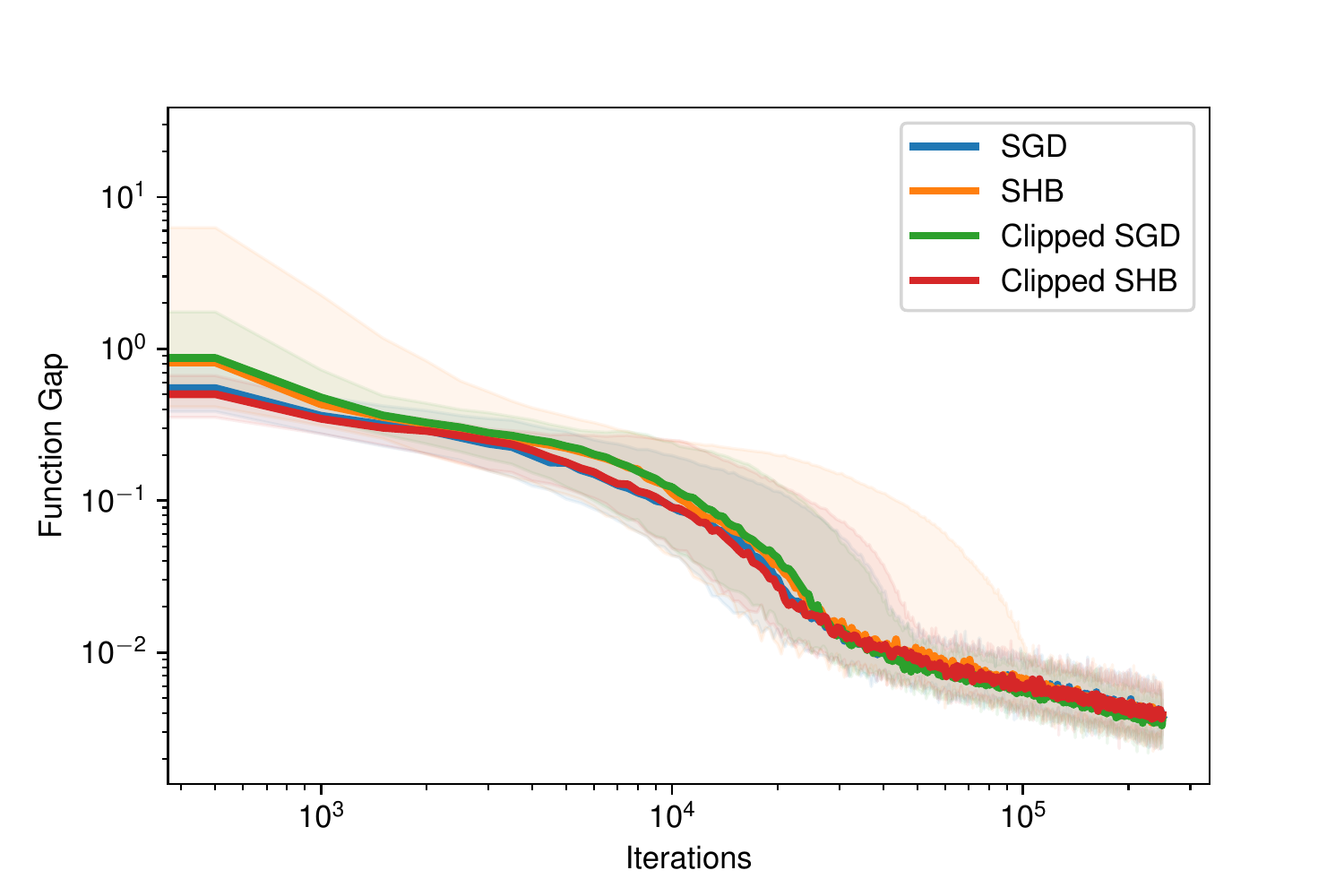}}
		\subcaption{$\stepsize_0 = 0.139$}	
		\end{minipage}
%	\hskip -0.15in
	\begin{minipage}{0.495\textwidth}
		\centering 
		{\includegraphics[width=1.\textwidth]{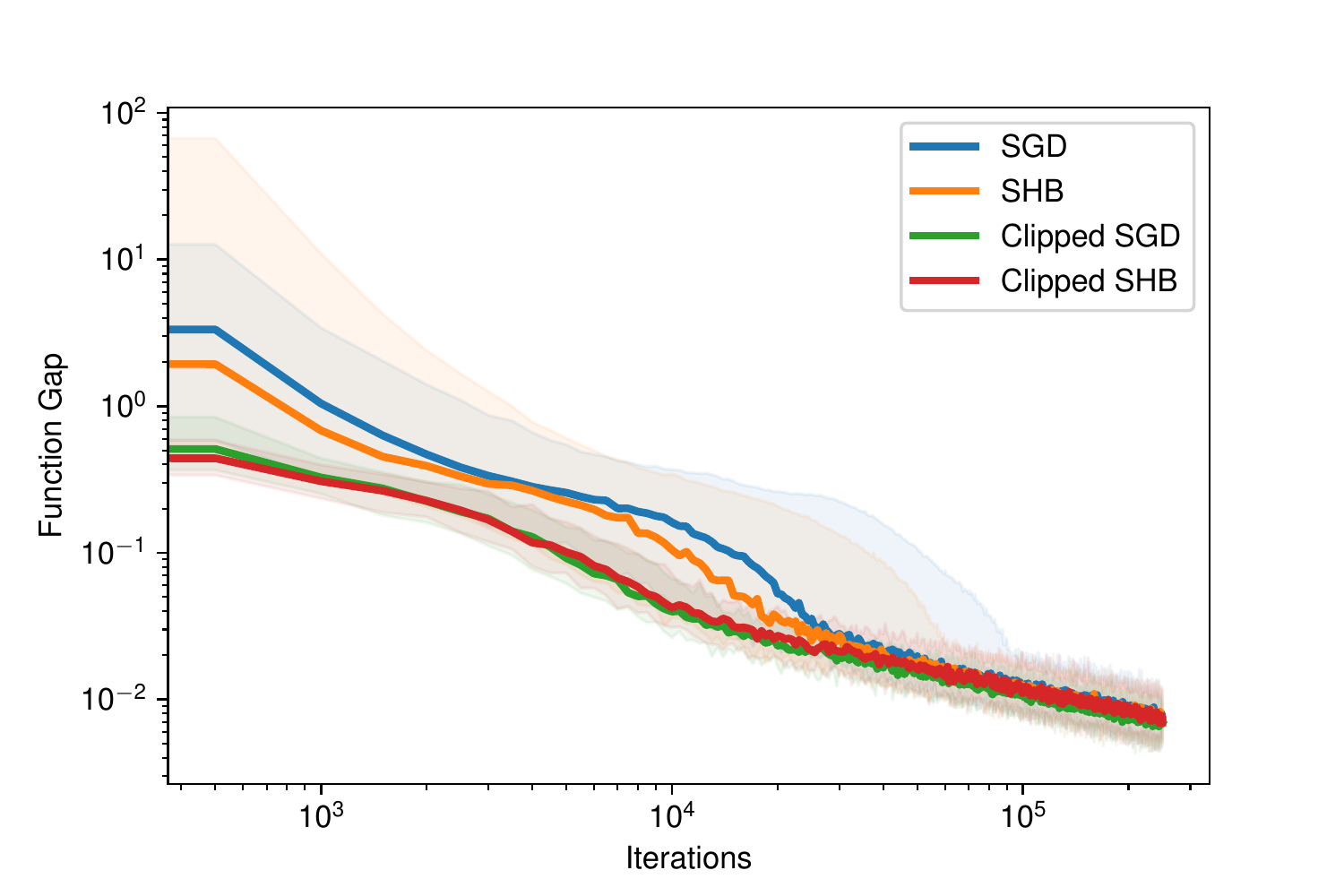}}
		\subcaption{$\stepsize_0=0.268$}	
	\end{minipage}
	%\hskip-0.15in
	\begin{minipage}{0.495\textwidth}
	\centering
	{\includegraphics[width=1.\textwidth]{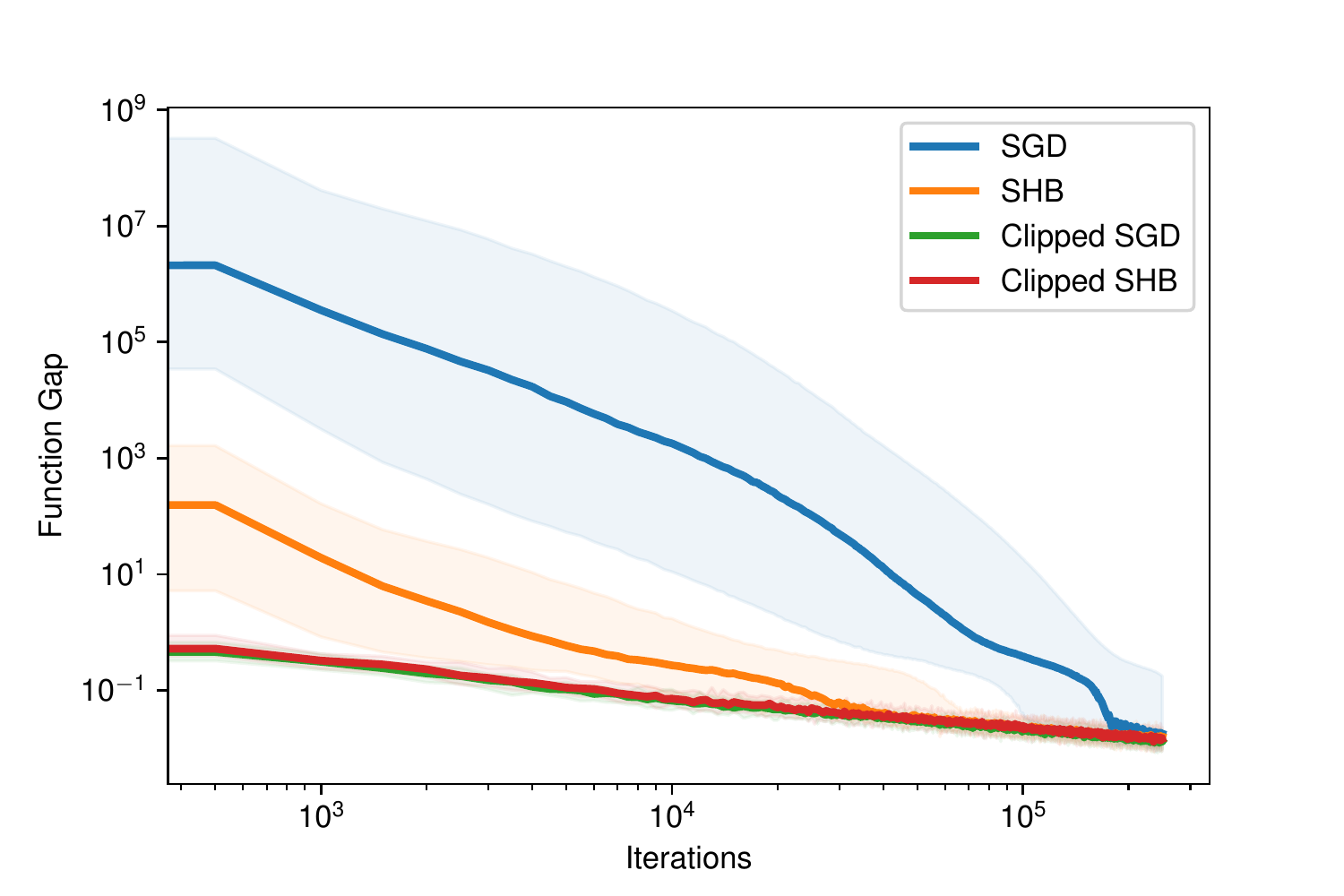}}
	\subcaption{$\stepsize_0=0.518$}	\label{fig:1:d}
	\end{minipage}
%	\hskip -0.15in
	\begin{minipage}{0.495\textwidth}
		\centering
		{\includegraphics[width=1.\textwidth]{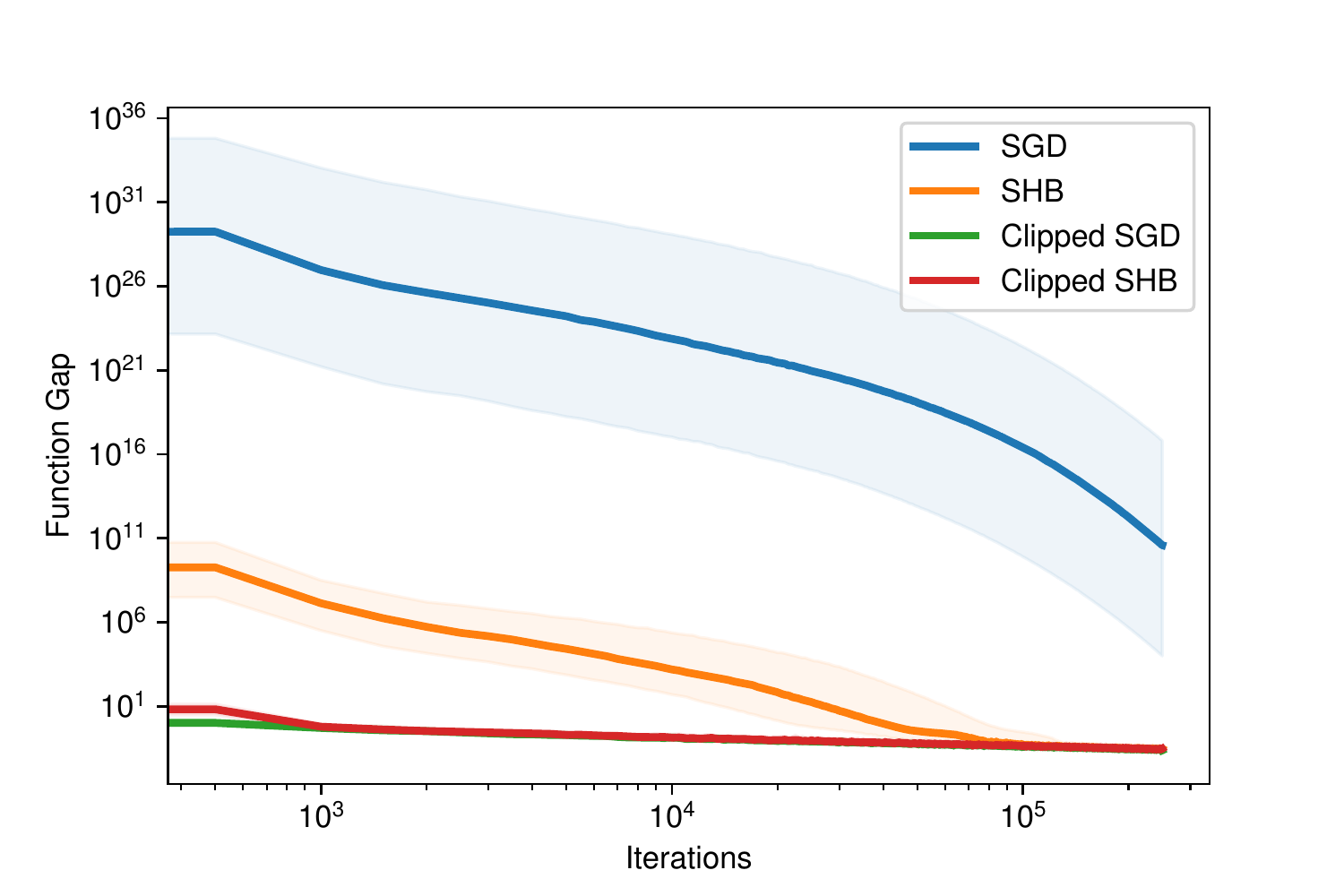}}
		\subcaption{$\stepsize_0=1.0$}	
	\end{minipage}	
		%\vskip -0.1cm
	\caption{The function gap $f(x_k)-f(x\opt)$ versus iteration count for phase retrieval with  $\gamma=10$.}\label{fig:PR}
\end{figure*}

%\vspace{-0.15cm}
\subsection{Phase retrieval}
Given $m$ measurements $(a_i, b_i) \in \R^{n}\times \R$, the (robust) phase retrieval problem seeks a vector $x\opt$ such that $\InP{a_i}{x\opt}^2 \approx b_i$ for most measurements  $i=1,\ldots,m$ by solving
\begin{align*}
	\minimize_{x\in\R^n} \frac{1}{m} \sum_{i=1}^{m} \big|\InP{a_i}{x}^2-b_i \big|.
\end{align*}
As for the vector $b$, in each problem instance, we select $x\opt$ uniformly from the unit sphere and construct its elements $b_i$ as $b_i = \InP{a_i}{x\opt}^2 + \delta \zeta_i$,  $i =1, \ldots,m,$ where $\zeta_i \sim \mc{N}(0, 25)$ models \replaced{corrupted measurements}{the corruptions}, and $\delta\in\{0,1\}$ is a binary random variable taking the value $1$ with probability $p_{\textrm{fail}}=0.1$, so that $p_{\textrm{fail}} \cdot m$ measurements are noisy. 

\begin{figure*}[h!]
	\centering
	\hskip-0.15in
	\begin{minipage}{0.35\textwidth}
		\centering
		{\includegraphics[width=1.\textwidth]{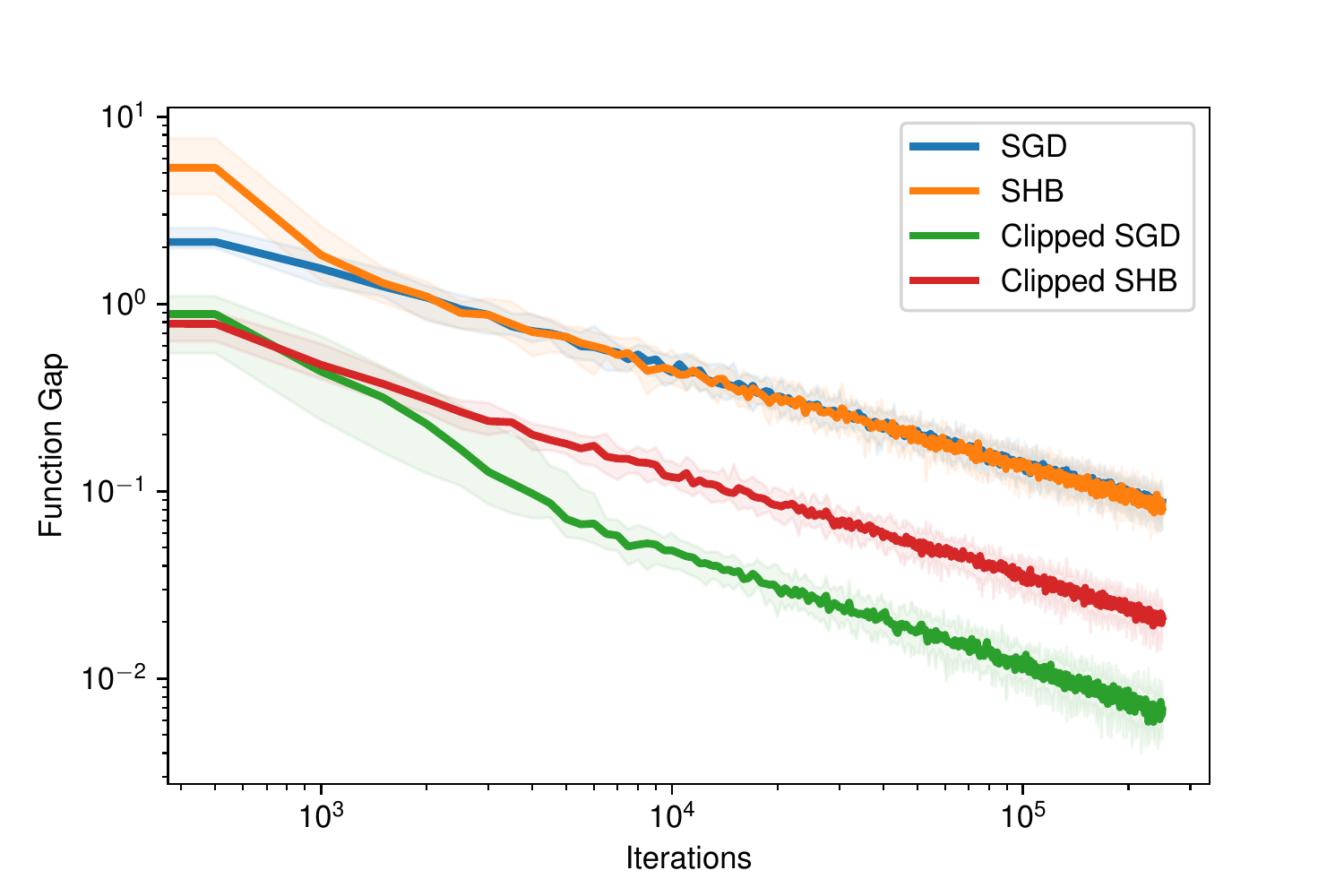}}
		\subcaption{$\gamma=0.5$}	
	\end{minipage}
	\hskip-0.15in
	\begin{minipage}{0.35\textwidth}
		\centering		{\includegraphics[width=1.\textwidth]{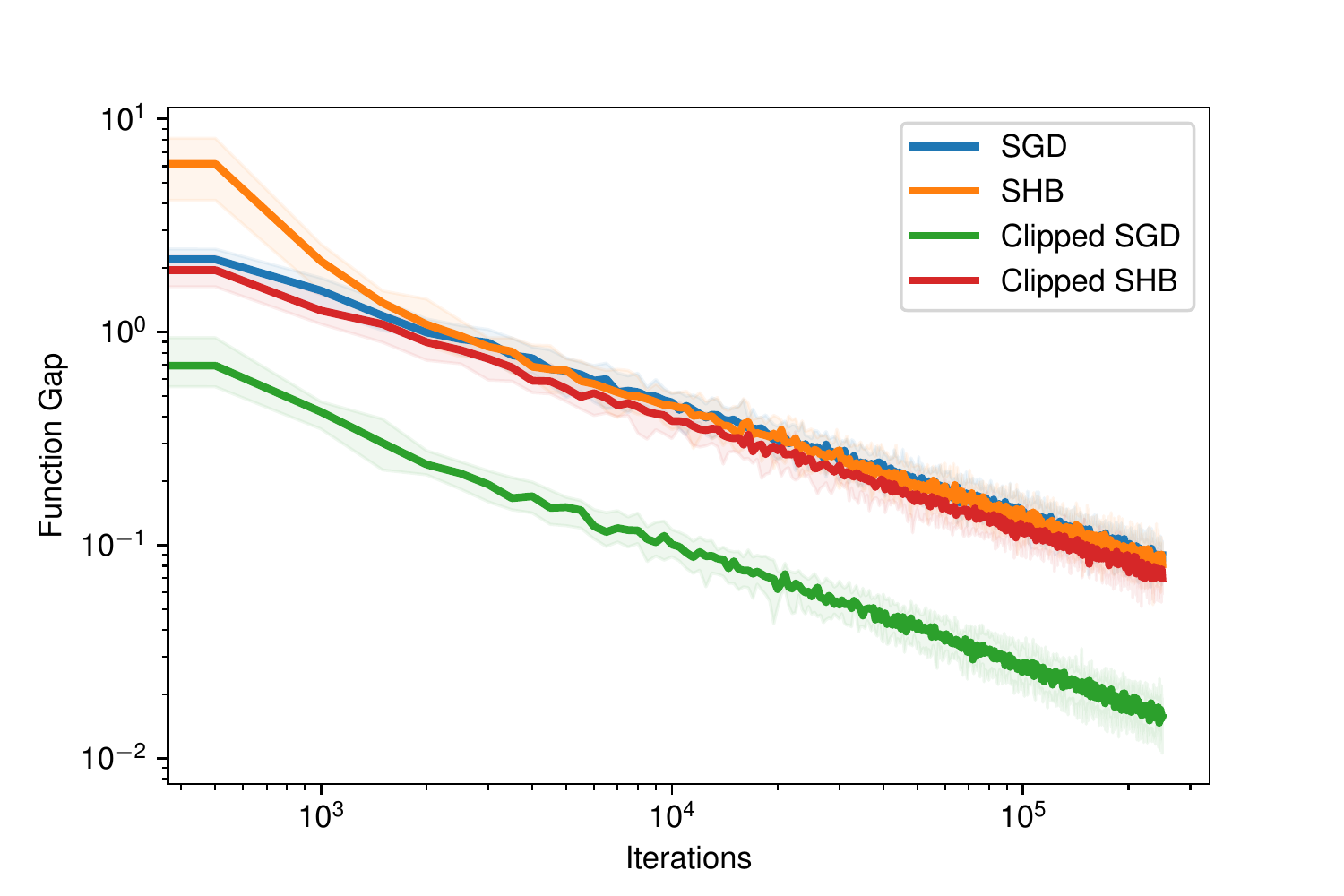}}
		\subcaption{$\gamma=1.0$}	
	\end{minipage}
	\hskip-0.15in
	\begin{minipage}{0.35\textwidth}
	\centering
	{\includegraphics[width=1.\textwidth]{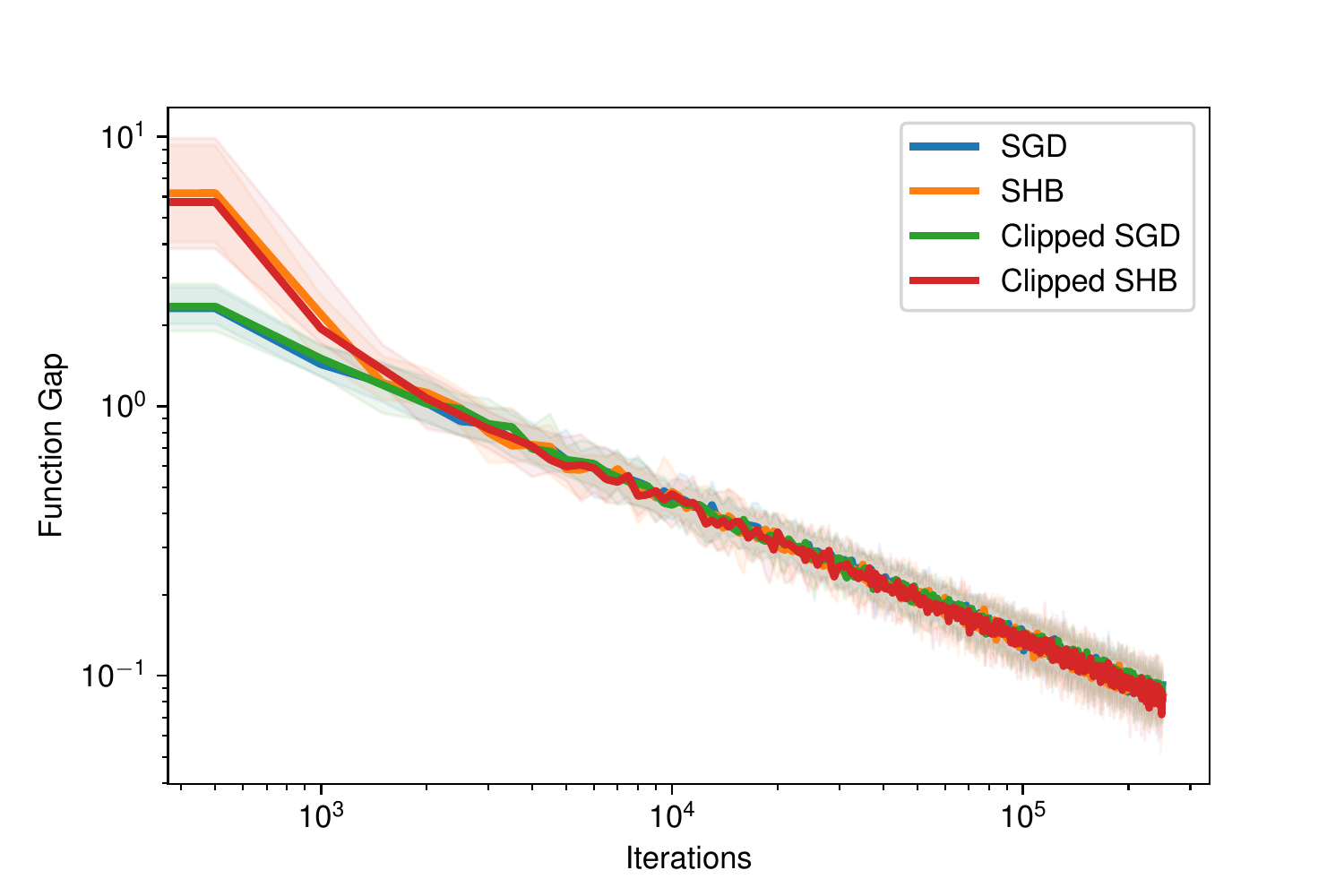}}
	\subcaption{$\gamma=10.0$}	
	\end{minipage}
		%\vskip -0.1cm
	\caption{The function gap $f(x_k)-f(x\opt)$ versus iteration count for absolute linear regression with $\stepsize_0=5$ and $1-\mmt=0.9$.}\label{fig:RLR}
	%\vspace{-0.2cm}
\end{figure*}
\begin{figure*}[h!]
%		\hskip -0.15in
	\centering
	\begin{minipage}{0.495\textwidth}
		{\includegraphics[width=1.\textwidth]{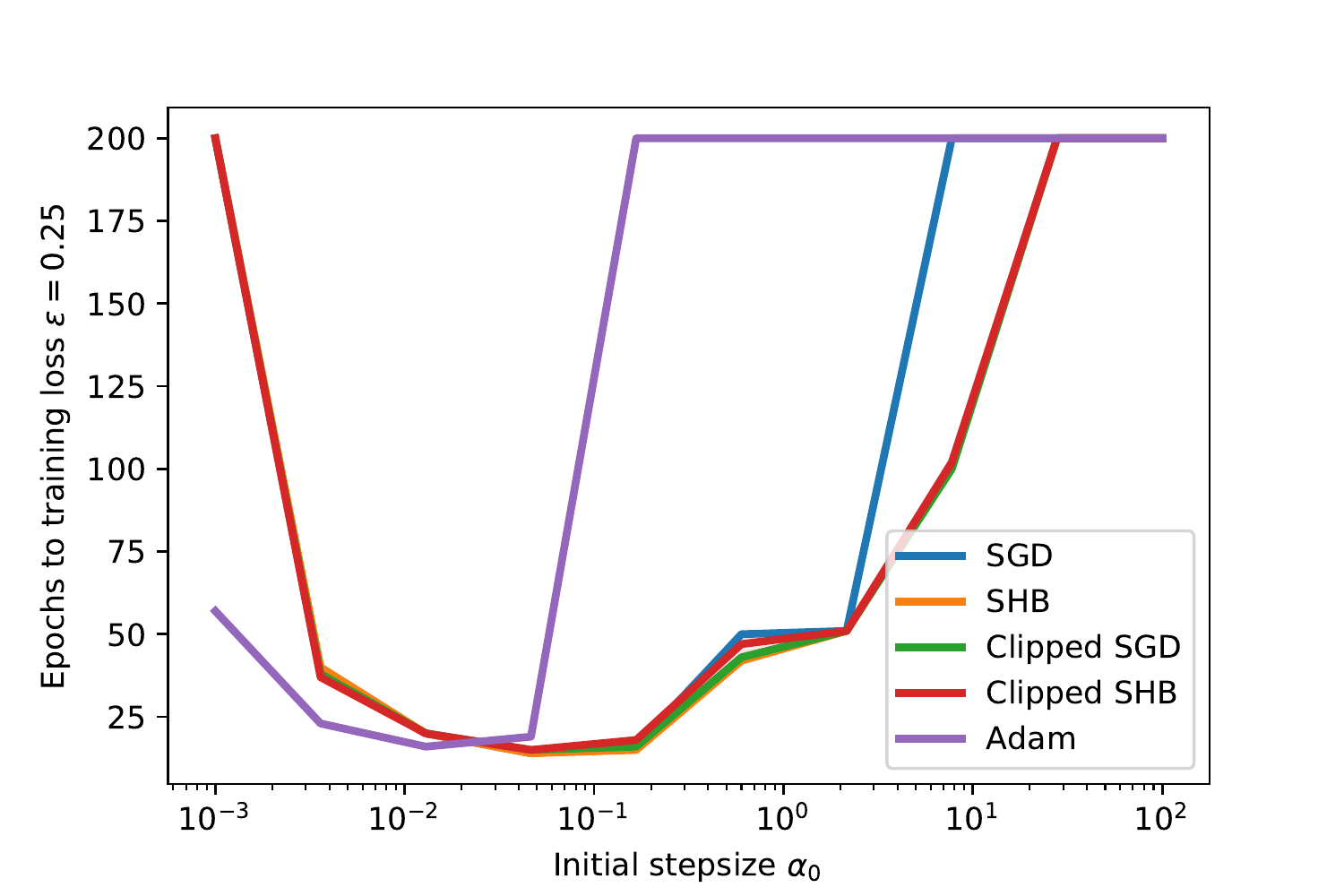}}
	\end{minipage}
%	\hskip-0.15in
	\begin{minipage}{0.495\textwidth}
		\centering
		{\includegraphics[width=1.\textwidth]{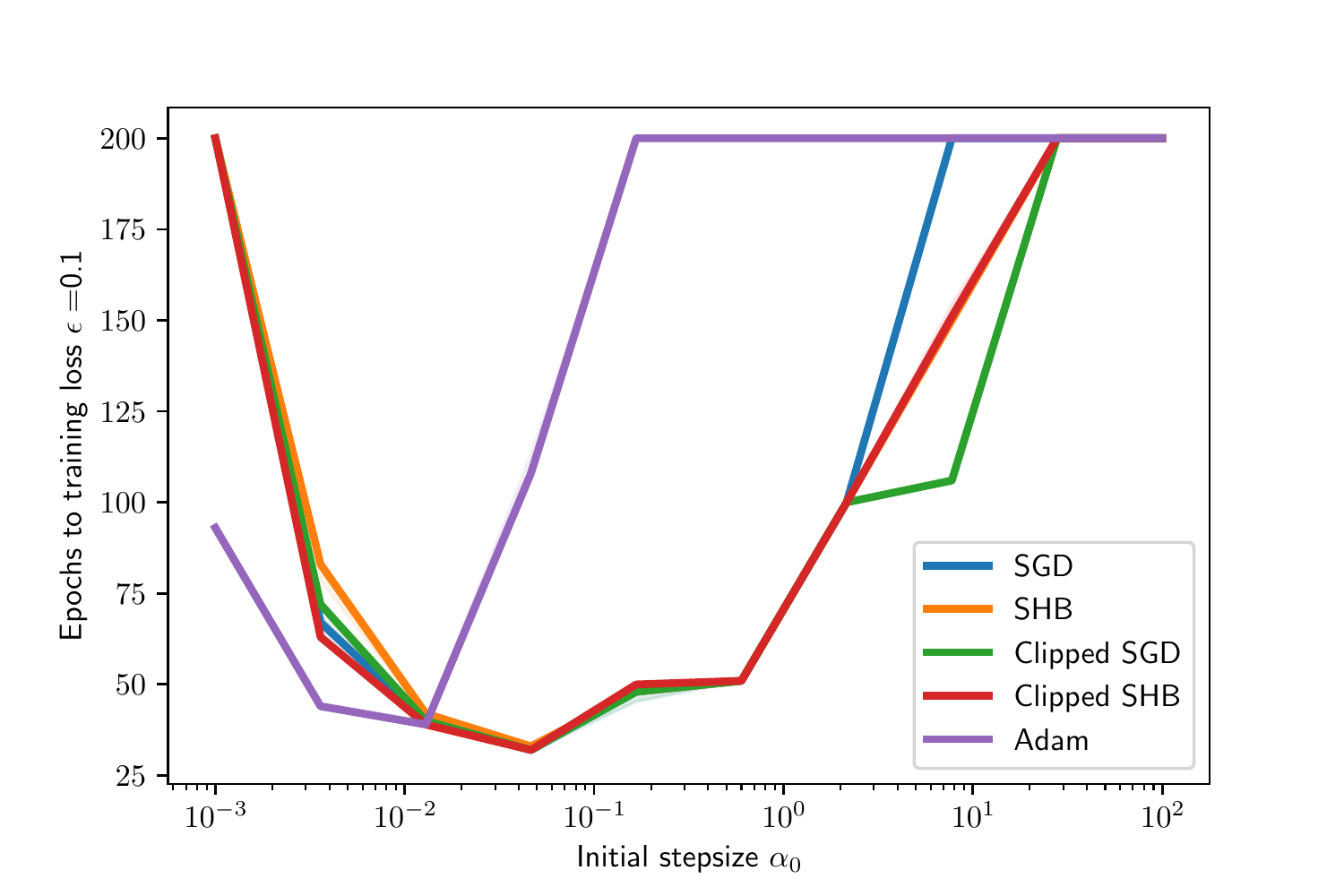}}
	\end{minipage}
	%\hskip-0.15in
%	\begin{minipage}{0.495\textwidth}
%		\centering
%		{\includegraphics[width=1.\textwidth]{{Figures/train_loss_clipthrs10.0_eps0.01}.pdf}}
%	\end{minipage}	
%	%\vskip -0.01in
%		%\hskip-0.15in			
%		\centering
	\begin{minipage}{0.495\textwidth}
		{\includegraphics[width=1.\textwidth]{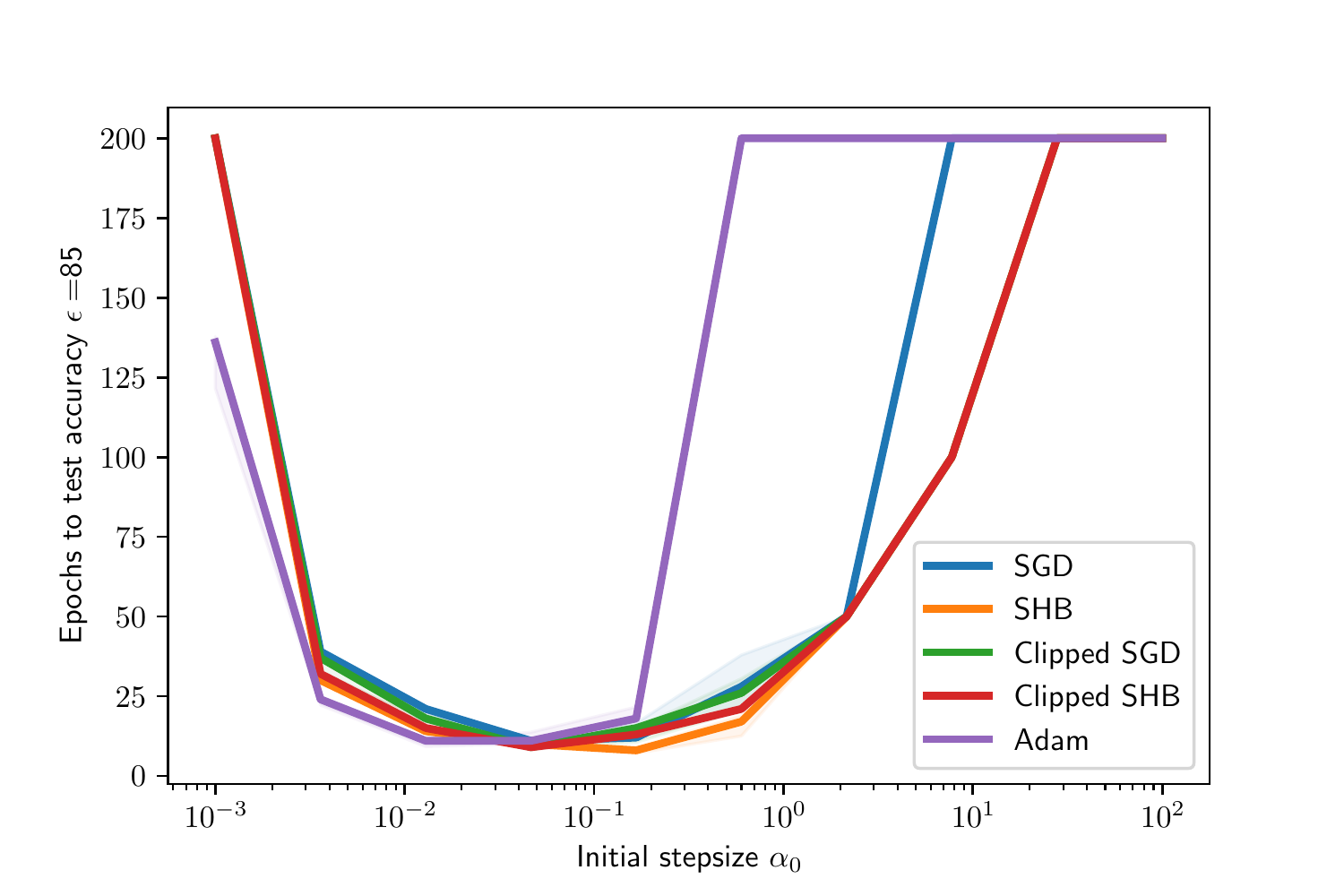}}
	\end{minipage}
	%\hskip-0.15in
	\begin{minipage}{0.495\textwidth}
		\centering	
		{\includegraphics[width=1.\textwidth]{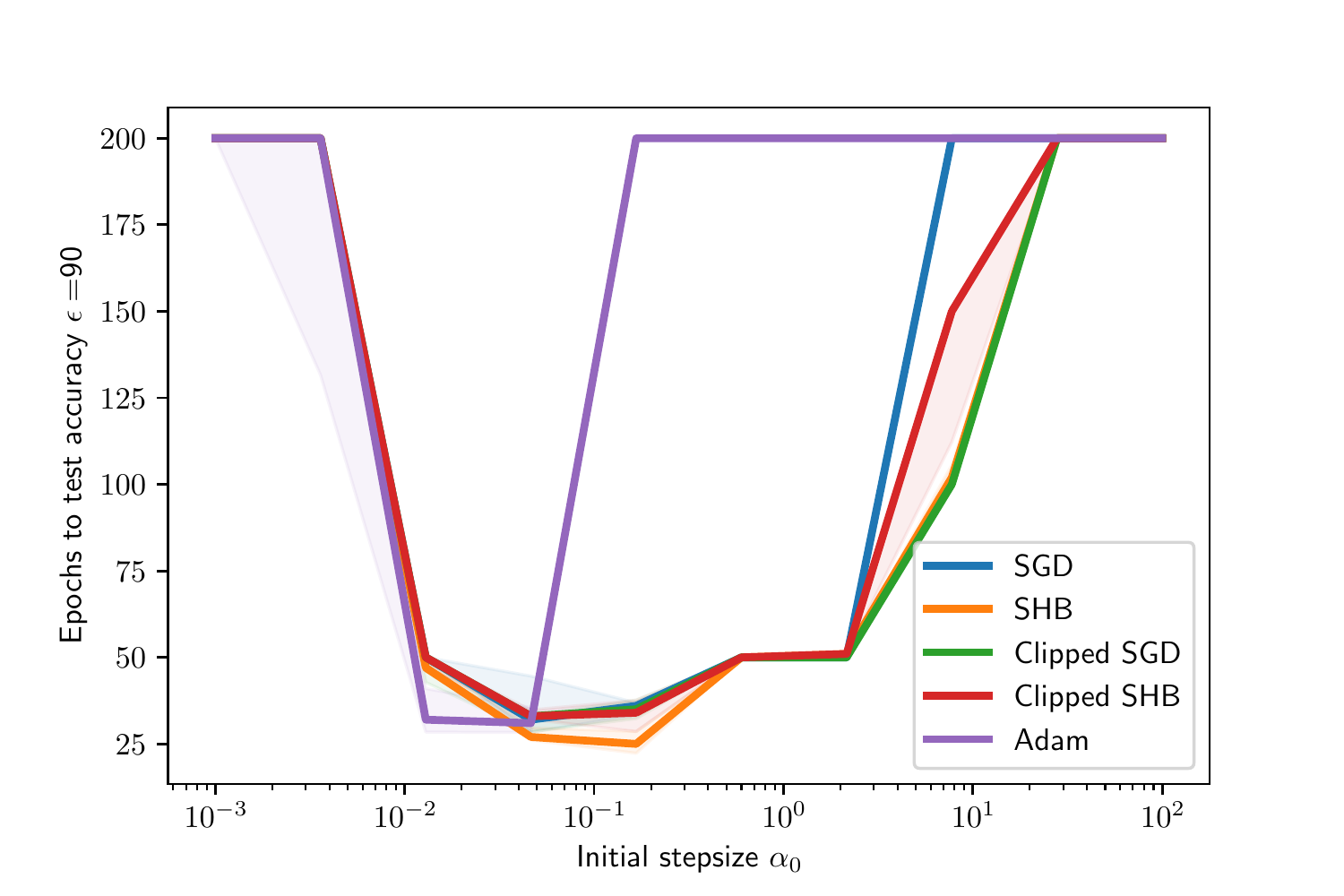}}
	\end{minipage}		
%	%\hskip-0.15in
%	\begin{minipage}{0.495\textwidth}
%		\centering 
%		{\includegraphics[width=1.\textwidth]{{Figures/test_acc_clipthrs10.0_eps93}.pdf}}
%	\end{minipage}	
	%\vskip -0.1cm
	\caption{The number of epochs to achieve $\epsilon$ training loss and test error versus initial stepsize $\stepsize_0$ for CIFAR10 with $\gamma=10$.}\label{fig:cifar:epoch2eps}
	%\vspace{-0.2cm}
\end{figure*}
\begin{figure}[h!]
	%\vskip -.25cm
	\centering
	{\includegraphics[width=0.7\textwidth]{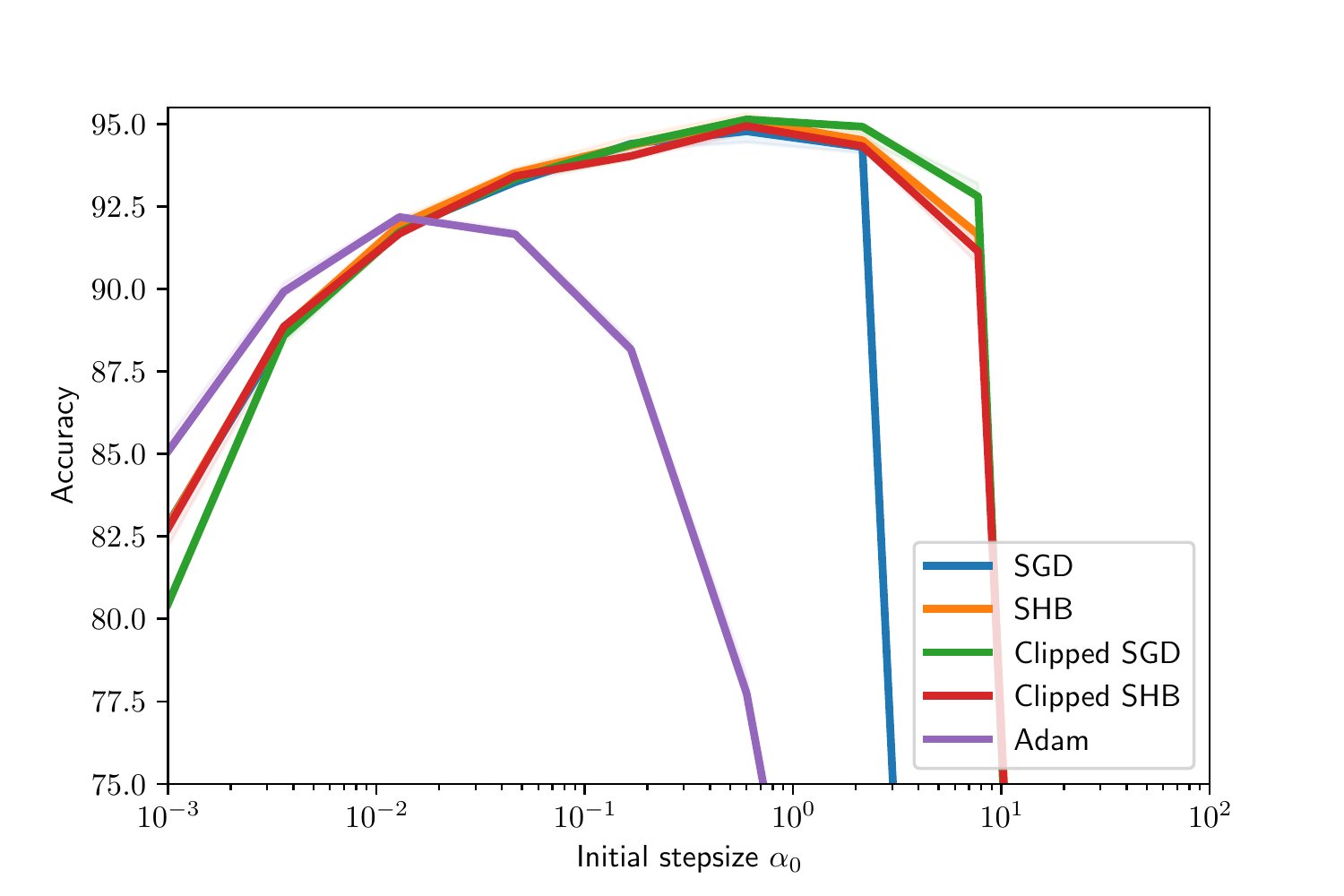}}
	%\vskip -0.1cm
\caption{The best achievable accuracy versus initial stepsize $\stepsize_0$ for CIFAR10 with $\gamma=10$.}\label{fig:cifar:best:acc}
%\vspace{-0.25cm}
\end{figure}
Figure~\ref{fig:PR:epoch2eps} shows the improved robustness provided by gradient clipping for the SGD and SHB algorithms. These clipped methods achieve good accuracies (within the allowed number of epochs) for much wider ranges of initial stepsizes than their unclipped versions. To further elaborate on this,  Figure~\ref{fig:PR} depicts the actual performance for 4 consecutive stepsizes (out of 15) used to produce Figure~\ref{fig:PR:epoch2eps}. We can see that these clipped methods always remain stable, while (started from the same initial point) SGD and SHB exhibit the problem of unboundedness when moving beyond their narrow ranges of working parameters.

%\vspace{-0.15cm}
\subsection{Absolute linear regression}
We consider mean absolute error $f(x) = \frac{1}{m}\lone{Ax-b}$. For each  problem instance, we generate $b=Ax\opt + \sigma w$ for $w\sim \mc{N}(0, 1)$ and $\sigma=0.01$. The problem is well-behaved as $f$ is both convex and Lipschitz continuous; we did not observe instability of SGD and SHB. Figure~\ref{fig:RLR} shows that gradient clipping does not harm  and sometimes can significantly boost the performance of their unclipped counterparts. We can see that although all methods converge with a similar slope, clipped methods may achieve better final accuracies. One possible explanation for this result would be the connection between the used stepsizes and the final error of the (stochastic) subgradient method; we have $\E[f(\bar{x}_k)]- f(x\opt) \leq O(\stepsize_k)$ for $\stepsize_k = O(1/\sqrt{k})$ \cite{BLA03, Duc18}.  With clipping, using a smaller $\gamma$ (if permitted) might have some effect on reducing the effective stepsizes $\stepsize_k \min\{1, \gamma/\ltwo{g_k}\}$, thereby yielding smaller errors.

%\vspace{-0.15cm}
\subsection{Neural Networks}
For our last set of experiments, we consider the image classification task on the CIFAR10 dataset \cite{KG09} with the ResNet-18 architecture \cite{HZRS16}. Here, we also compare the previous methods with the Adam algorithm using its default parameters in PyTorch; $\beta_1=0.9$, $\beta_2=0.99$, and $\epsilon=10^{-8}$.\footnote{\url{https://pytorch.org}} 
Following common practice, we use mini-batch size 128, momentum parameter $\beta=0.9$, and weight-decay coefficient $5\times 10^{-4}$ in all experiments.  \replaced{For each algorithm, we conduct}{We conduct for each algorithm} 5 experiments (up to 200 epochs) and report the medians of the training loss and test accuracy together with the 90\% confidence intervals. For the stepsizes, we use constant values starting with $\stepsize_0$ and reduce them by a factor of $10$ every 50 epochs. The initial stepsizes $\stepsize_0$ for Adam are scaled by $1/100$ in actual runs \cite{AD19b}.

Figure~\ref{fig:cifar:epoch2eps} shows the minimum number of epochs required to reach desired values for various performance measures as a function of the initial stepsize. As the classification task on CIFAR10 is a rather well-conditioned problem, \replaced{the results tell a very similar story to our absolute linear regression experiments.}{we have a very similar story to the absolute linear regression one.} We also observe that Adam is more sensitive to stepsize selection and needs more time to achieve good test performance in this example.  To further clarify this, Figure~\ref{fig:cifar:best:acc} shows that over the tested range of stepsizes, Adam is not able to reach the same best achievable test accuracies that the other methods do.

In summary, the results in this section reinforce our theoretical developments that gradient clipping can: i) stabilize and guarantee convergence for problems with rapidly growing gradients; ii) retain and sometimes improve the best performance of their unclipped counterparts even on standard (``easy'')  problems. 

\section{Conclusions}
We analyzed clipped subgradient-based methods for solving stochastic convex and non-convex optimization problems. Moving beyond traditional quadratic models, we showed that these methods enjoy strong stability properties and attain classical convergence rates in settings where standard convergence theory does not apply. With a novel Lyapunov analysis, we also proved that the sample complexity of the methods match  the best-known result for weakly convex problems, emphasizing the effectiveness of gradient clipping on a wide range of problem classes.

\section*{Acknowledgement}
This work was supported in part by the Knut and Alice Wallenberg Foundation, the Swedish Research Council and the Swedish Foundation for Strategic Research. 
The computations were enabled by resources provided by the Swedish National Infrastructure for Computing (SNIC), partially funded by the Swedish Research Council through grant agreement no. 2018-05973. Finally,  we thank the anonymous reviewers for their detailed and valuable feedback.

\bibliographystyle{abbrv}
\bibliography{refs}

\newpage
\appendix

\noindent In the sequel, we will frequently use the following Young's inequality
\begin{align}\label{eq:young}
	\abs{\InP{a}{b}}\leq \frac{\ltwo{a}^2}{2\epsilon} + \frac{\epsilon \ltwo{b}^2}{2},
\end{align}
which holds for any $a,b\in\R^n$ and $\epsilon>0$.

\section{Proof of Example~\ref{example:sgd:divergence}}\label{appendix:proof:example:sgd:divergence}
Recall that $\stepsize_k=\stepsize_1/k$ and
\begin{align*}
	x_{k+1} = x_k - \frac{\stepsize_1}{k} \left(x_k^3 + \epsilon x_k \right).
\end{align*}
We follow \cite{And90} and prove our claim by induction. The hypothesis is obviously true for  $k=1$. For $k=2$, since $x_1\geq \sqrt{3/\stepsize_1}$, it holds that
\begin{align*}
	\abs{x_2} = \abs{x_1} \abs{1 - \stepsize_1 \left(x_1^2 + \epsilon\right)} = \abs{x_1} \left(\stepsize_1 \left(x_1^2 + \epsilon\right) -1\right) \geq 2 \abs{x_1}.
\end{align*} 
Suppose that the hypothesis is true up to some iteration $k>2$, we have
\begin{align*}
	\abs{x_{k+1}} 
	&= \abs{x_k} \abs{1 - \frac{\stepsize_1}{k} \left(x_{k}^2 + \epsilon\right)} 
	= \abs{x_k} \left(\frac{\stepsize_1}{k} \left(x_{k}^2 + \epsilon\right) -1\right) 
	\geq  \abs{x_k} \left(\frac{\stepsize_1}{k} x_{k}^2 -1\right) 
	\nonumber\\
	&\geq  \abs{x_k} \left(\frac{3}{k} (k!)^2 -1\right) 
	\geq \abs{x_k} \left(3k(k-1)- 1\right)
	\geq \abs{x_k} (k+1) 
	\geq \abs{x_1} (k+1)!,
\end{align*} 
where we have used the facts that $(k!)^2/k = k! (k-1)! \geq k(k-1) $ and $3k(k-1)- 1 \geq k+1$ for any $k\geq 2$.

\section{Proof of Proposition~\ref{propositiom:stability}}\label{appendix:proof:proposition:stability}
To obtain \eqref{eq:proposition:stability}, it is sufficient to show that
\begin{align*}%\label{eq:almost:quasi:momotone}
	\E\left[ \dist\left(x_{k+1}, \Xopt\right)^2 \big | \mc{F}_{k} \right]  \leq \dist\left(x_{k}, \Xopt\right)^2 + \stepsize_k C
\end{align*}
for some constant $C$.  Let $x\opt = \gproj{\Xopt}{x_k}$, we have
\begin{align} \label{eq:proposition:stability:proof:xdiff:1}
	\ltwo{x_{k+1}-x\opt} ^2
	= 
		 \ltwo{x_k-x\opt} ^2
		- 2 \stepsize_k \InP{d_k}{x_k-x\opt}
		+ \stepsize_k^2 \ltwo{d_k}^2.
\end{align}
Let $\varrho_k := \min\left\{1, \gamma_k/\ltwo{g_k}\right\} $, so that $d_k=\varrho_k g_k$.
The established approaches develop the preceding inner product by adding and subtracting either the true subgradient $f'(x_k)$ or its clipped version, which is  likely to require an upper bound for $\ltwo{f'(x_k)}$. In contrast, we add and subtract the term $\varrho_k f'(x_k)$ to obtain 
\begin{align*}
	- \stepsize_k \InP{d_k}{x_k-x\opt}
	=	
		\stepsize_k \varrho_k \InP{f'(x_k) - g_k} {x_k-x\opt}
		-\stepsize_k  \varrho_k \InP{f'(x_k)}{x_k-x\opt}.
\end{align*} 
Now, by the convexity of $f$ and Assumption~\ref{assumption:quadratic:growth}, we have 
\begin{align*}
 \InP{f'(x_k)}{x_k-x\opt} \geq f(x_k) - f(x\opt) \geq \mu \ltwo{x_k-x\opt}^2.
\end{align*}
It follows that
\begin{align}\label{eq:proposition:stability:proof:inp}
	- \stepsize_k \InP{d_k}{x_k-x\opt}
	&\leq	
		\stepsize_k \varrho_k 
		\left(	\InP{\grad{f(x_k)} - g_k }{x_k-x\opt}
				- \mu \ltwo{x_k-x\opt}^2
		\right)
	\nonumber\\
	&\leq
		\frac{\stepsize_k \varrho_k}{4\mu} \ltwo{f'(x_k) - g_k }^2
	\nonumber\\
	&\leq	\frac{\stepsize_k }{4\mu} \ltwo{f'(x_k) - g_k }^2,
\end{align} 
where we used Young's inequality \eqref{eq:young} with $\epsilon=2\mu$, and the fact that $\varrho_k\leq 1$.
Combining \eqref{eq:proposition:stability:proof:inp} and \eqref{eq:proposition:stability:proof:xdiff:1} and noticing that $\stepsize_k^2 \ltwo{d_k}^2 \leq \stepsize_k^2\gamma_k^2 \leq \gamma^2 \stepsize_k$, we obtain
\begin{align} \label{eq:proposition:stability:proof:xdiff:2}
	\ltwo{x_{k+1}-x\opt} ^2
	\leq
		 \ltwo{x_k-x\opt} ^2
		+\frac{\stepsize_k }{2\mu} \ltwo{f'(x_k) - g_k }^2
		+ \gamma^2\stepsize_k .
\end{align}
Taking the conditional expectation in \eqref{eq:proposition:stability:proof:xdiff:2} and using Assumption~\ref{assumption:variance:cvx}  yields
%\begin{align}\label{eq:proposition:stability:proof:inp:2}
%-\stepsize_k \E\left[\InP{d_k}{x_k-x\opt} \big | \mc{F}_k\right]
%	\leq	
%		\frac{\stepsize_k }{4\mu}
%		\E\left[ \ltwo{f'(x_k) - g_k }^2 \big | \mc{F}_k\right]
%	\leq 
%			\frac{\sigma^2 }{4\mu} \stepsize_k.
%\end{align} 
%Since $\gamma_k \leq \gamma / \stepsize_k$, we have $\stepsize_k^2 \ltwo{d_k}^2 \leq  \gamma^2 \stepsize_k$. Now, plugging this and \eqref{eq:proposition:stability:proof:inp:2} into \eqref{eq:proposition:stability:proof:xdiff:1}, we arrive at
\begin{align*}
	\E\left[ \ltwo{x_{k+1}-x\opt} ^2 \big | \mc{F}_k \right]
	= 
		 \ltwo{x_k-x\opt} ^2
		+ 
		\left({\sigma^2}/(2\mu) + \gamma^2 \right) \stepsize_k,
\end{align*}
which completes the proof since $\dist\left(x_{k+1}, \Xopt\right) \leq \ltwo{x_{k+1}-x\opt} ^2$ for every $x\opt\in \Xopt$.

\section{Proof of Theorem~\ref{thrm:cvx:minibatch:asym}}\label{appendix:proof:thrm:cvx:minibatch:asym}

In this section, we prove Theorem~\ref{thrm:cvx:minibatch:asym}. We begin with the following classical result  on the almost sure convergence of nonnegative almost supermatingales.
\begin{lemma}[Robbins-Siegmund \cite{RS71}]\label{lem:RS}
Let $A_k, B_k, C_k$, and $ V_k$ be non-negative random variables adapted to the filtration $\mc{F}_k$ and satisfying:
\begin{align*}
	\E\left[V_{k+1}|\mc{F}_k\right] \leq (1+A_k)V_k + B_k - C_k, \quad \forall k \in \N.
\end{align*}
Then, on the event $\left\{\sum_{k=0}^{\infty} A_k < \infty, \sum_{k=0}^{\infty} B_k < \infty \right\}$,
there is a random variable $V_\infty$ such that 
\begin{align*}
	V_k \lcas V_\infty <\infty
	\quad \mbox{and} \quad 
	\sum_{k=0}^{\infty} C_k < \infty \quad \mbox{a.s.}
\end{align*}
\end{lemma}
\medskip

\noindent We can now proceed as follows. Let $\varrho_k := \min\left\{1, \gamma/\ltwo{g_k}\right\} $ and $x\opt= \gproj{\Xopt}{x_k}$, it holds that
\begin{align*}
	 \ltwo{x_{k+1}-x\opt} ^2
	 =
			 \ltwo{x_k-x\opt} ^2
			- 2\stepsize_k \varrho_k  \InP{f'(x_k)}{x_k-x\opt}
			+ 2\stepsize_k \varrho_k  \InP{ f'(x_k) - g_k}{x_k-x\opt}
			+ \stepsize_k^2 \ltwo{d_k}^2.
\end{align*}
By the convexity of $f$, Assumption~\ref{assumption:quadratic:growth}, and Young's inequality \eqref{eq:young} with $\epsilon=\mu$, we have 
\begin{align*}
 \InP{f'(x_k)}{x_k-x\opt} &\geq f(x_k) - f(x\opt) \geq \mu \ltwo{x_k-x\opt}^2,
 \nonumber\\
 \InP{ f'(x_k) - g_k}{x_k-x\opt} &\leq \frac{1}{2\mu}\ltwo{f'(x_k) - g_k}^2 + \frac{\mu}{2}\ltwo{x_k -x\opt}^2.
\end{align*}
We thus arrive at
\begin{align*}
	\ltwo{x_{k+1}-x\opt} ^2
	&\leq 
		\left(1 - \mu \stepsize_k\varrho_k  \right)
	 \ltwo{x_k-x\opt} ^2
		+ \frac{\stepsize_k\varrho_k }{\mu}\ltwo{f'(x_k) - g_k}^2
	+ \stepsize_k^2 \gamma^2,			
%	\nonumber\\
%	&\leq 
%	\left(1 - \mu \stepsize_k\varrho_k \right)
%	 \ltwo{x_k-x\opt} ^2
%		+ \frac{\stepsize_k}{\mu}\ltwo{f'(x_k) - g_k}^2
%		+ \stepsize_k^2 \gamma^2,
\end{align*}
where we have used $\ltwo{d_k}\leq \gamma$. Since Assumption~\ref{assumption:variance:cvx} and the definition of $g_k$ implies  $\E\left[\ltwo{f'(x_k) - g_k}^2 \big | \mc{F}_{k} \right] \leq \sigma^2/\batchsize_k$,
taking the conditional expectation in the preceding inequality thus yields
\begin{align*}
		\E\left[ \dist\left(x_{k+1}, \Xopt\right)^2 \big | \mc{F}_{k} \right] 
%		\leq	
%			\E\left[\ltwo{x_{k+1}-x\opt} ^2 \big | \mc{F}_k \right]
% 		&\leq 
%			\left(1 - \mu \stepsize_k 	\E\left[\varrho_k \big | \mc{F}_k \right]  \right)
%			 \ltwo{x_k-x\opt} ^2
%			+ \frac{\stepsize_k}{\mu}	\E\left[\ltwo{f'(x_k) - g_k}^2 \big | \mc{F}_k \right]
%			+ \stepsize_k^2 \gamma^2
%		\nonumber\\
		\leq
			\left(1 - \mu \stepsize_k 	\E\left[\varrho_k \big | \mc{F}_k \right] \right)
			 \dist\left(x_{k}, \Xopt\right)^2 
			+ \frac{ \sigma^2 \stepsize_k}{\mu \batchsize_k} 
			+ \stepsize_k^2 \gamma^2.
\end{align*}
Now, invoking Lemma~\ref{lem:RS} with $V_k =  \dist\left(x_{k}, \Xopt\right)^2$, $A_k= 0$, $B_k=\frac{ \sigma^2 \stepsize_k}{\mu \batchsize_k} + \stepsize_k^2 \gamma^2$, and $C_k=\mu \stepsize_k 	\E\left[\varrho_k \big | \mc{F}_k \right]  \dist\left(x_{k}, \Xopt\right)^2$, we have 
\begin{align*}
	V_k \lcas V_\infty < \infty 
	\quad \mbox{and} \quad
	\sum_{k=0}^{\infty} \stepsize_k 	\E\left[\varrho_k \big | \mc{F}_k \right]  \dist\left(x_{k}, \Xopt\right)^2 < \infty \quad \mbox{a.s.}
\end{align*}
Next, we will show that $\E\left[\varrho_k \big | \mc{F}_k \right]$ is bounded away from zero almost surely.
Since $\min\{a,b\} = ab/\max\{a,b\} \geq ab/(a+b)$ for  $a,b\geq 0$, it holds that
\begin{align*}
	\E\left[\varrho_k \big | \mc{F}_k \right] 
	\geq
		\E\left[ \frac{\gamma}{\gamma + \ltwo{g_k}}  \bigg | \mc{F}_k \right] 
	\mathop \geq \limits^{\mathrm{(a)}} 
		 \frac{\gamma}{\gamma + \E\left[\ltwo{g_k} \big | \mc{F}_k \right] }
	\mathop \geq \limits^{\mathrm{(b)}}  
		\frac{\gamma}{\gamma + \left(\E\left[\ltwo{g_k}^2 \big | \mc{F}_k \right]\right)^{1/2} },
\end{align*}
where $\mathrm{(a)}$ and $\mathrm{(b)}$ follow from Jensen's inequality applied to the convex  function $1/(1+x)$ for $x\geq 0$ and the concave function $x^{1/2}$, respectively.  Thus, by Assumption~\ref{assumption:Gbig} and the definition of $g_k$, we have
\begin{align}\label{hrm:cvx:minibatch:asym:varro}
	\E\left[\varrho_k \big | \mc{F}_k \right] 
	\geq
		 \frac{\gamma}{\gamma + G_{\mathrm{big}}^{1/2}(\dist(x_k, \Xopt))}.
\end{align}
Now, since $V_k \lcas V_\infty < \infty$, $\dist(x_k, \Xopt)$ is bounded with probability one. It follows that $\E\left[\varrho_k \big | \mc{F}_k \right] $ is bounded away from zero almost surely. Consequently, $\sum_{k=0}^{\infty} \stepsize_k \dist\left(x_{k}, \Xopt\right)^2 < \infty$ almost surely. Since  $\sum_{k=0}^{\infty} \stepsize_k=\infty$, we deduce that $\liminf_{k\to \infty} \dist\left(x_{k}, \Xopt\right) = 0$ almost surely. But we know that the limit exists (since $V_k \lcas V_\infty < \infty$), and hence it must be the case that $ \dist\left(x_{k}, \Xopt\right) \lcas 0$.

\section{Proof of Theorem~\ref{thrm:cvx:minibatch:finite}}\label{appendix:proof:thrm:cvx:minibatch:finite}

We start by defining the following quantity for some real constant $M>0$:
\begin{align*}
	s = \inf\left\{k \geq 0 : \dist(x_k, \Xopt) > M\right\}.
\end{align*}
Let $e_k =  \dist\left(x_{k}, \Xopt\right)$, then by Theorem~\ref{thrm:cvx:minibatch:asym}, we have
\begin{align*}
	\E\left[e_{k+1}^2 \big | \mc{F}_{k} \right] 
		\leq
			\left(1 - \mu \stepsize_k 	\E\left[\varrho_k \big | \mc{F}_k \right] \right)
			  e_k^2
			+ \frac{ \sigma^2  \stepsize_k}{\mu \batchsize_k}
			+ \stepsize_k^2 \gamma^2.
\end{align*}
Defining the series of events $E_k$ as $E_k := \{s > k\}$, the preceding inequality implies
\begin{align}\label{eq:thrm:cvx:minibatch:finite:proof:e:1}
	\E\left[ e_{k+1}^2 \bindic{E_{k+1}} \big | \mc{F}_{k} \right] 
	&\leq 
		\E\left[ e_{k+1}^2 \bindic{E_{k}} \big | \mc{F}_{k} \right] 
	\nonumber\\
	&\leq
		\left(1 - \mu \stepsize_k 	\E\left[\varrho_k \big | \mc{F}_k \right] \right)
  e_k^2 \bindic{E_{k}}
		+ \frac{ \sigma^2  \stepsize_k}{\mu \batchsize_k}
		+ \stepsize_k^2 \gamma^2,	
\end{align}
where the first inequality follows since $E_{k+1} \subset E_k$, and the second one holds since $E_k \in \mc{F}_k$. We next bound the quantity of interest $\Pr\left(e_K^2 \leq \epsilon \right) $ as 
\begin{align}\label{eq:thrm:cvx:minibatch:finite:proof:Pe:1}
	\Pr\left(e_K^2 \leq \epsilon \right) 
	\geq 
		\Pr\left( e_K^2 \leq \epsilon\,\,\, \mbox{and}\,\,\, E_K\right)
%	\nonumber\\
	&= 
		\Pr\left( e_K^2 \leq \epsilon \, | E_K\right) \Pr\left(E_K\right)
	\nonumber\\
	&= 
		\left(1 - \Pr\left( e_K^2 > \epsilon \, | E_K\right)\right) \Pr\left(E_K\right)
	\nonumber\\
	&\geq
			\Pr\left(E_K\right) - \frac{\E\left[e_K^2 \, | E_K\right] }{\epsilon} \Pr\left(E_K\right)
	\nonumber\\
	&=
		\Pr\left(E_K\right) - \frac{\E\left[e_K^2 \bindic{E_K}\right] }{\epsilon} ,
\end{align} 
where we used the Markov's inequality and basic manipulations.
Next, we upper bound the term $\E\left[e_K^2 \bindic{E_K}\right]$. Since on the event $E_K$, $\dist\left(x_k, \Xopt\right) \leq M$ for any $k\leq K$, it follows from \eqref{hrm:cvx:minibatch:asym:varro} that
\begin{align}\label{eq:thrm:cvx:minibatch:finite:proof:varrho}
	\E\left[\varrho_k \big | \mc{F}_k \right] 
	\geq
		 \frac{\gamma}{\gamma + G_{\mathrm{big}}^{1/2}(\dist(x_k, \Xopt))}
	\geq 
		 \frac{\gamma}{\gamma + G_{\mathrm{big}}^{1/2}(M)}
	:= \varrho. 
\end{align}
In view of \eqref{eq:thrm:cvx:minibatch:finite:proof:e:1} and \eqref{eq:thrm:cvx:minibatch:finite:proof:varrho}, together with $m_k=1/\stepsize_k$, it holds that
\begin{align}\label{eq:thrm:cvx:minibatch:finite:proof:e:2}
	\E\left[ e_{k+1}^2 \bindic{E_{k+1}} \big | \mc{F}_{k} \right] 
	\leq
		\left(1 - \mu \varrho \stepsize_k  \right) 
		e_k^2 \bindic{E_{k}}
		+ \left({ \sigma^2  }/{\mu } +  \gamma^2\right)\stepsize_k^2.
\end{align}
By successively applying \eqref{eq:thrm:cvx:minibatch:finite:proof:e:2} and using $\stepsize_k=\stepsize_0(k+1)^{-\tau}$ with $\tau\in(1/2,1)$, we obtain
\begin{align}\label{eq:thrm:cvx:minibatch:finite:proof:e:3}
	\E\left[ e_{K}^2 \bindic{E_{K}}\right] \leq c_0 \stepsize_K,
\end{align}
where $c_0$ is some numerical constant \cite[Eq.~(A14)]{PJ92}.
Therefore, plugging \eqref{eq:thrm:cvx:minibatch:finite:proof:e:3} into \eqref{eq:thrm:cvx:minibatch:finite:proof:Pe:1} yields
\begin{align}\label{eq:thrm:cvx:minibatch:finite:proof:Pe:2}
	\Pr\left(e_K^2 \leq \epsilon \right) 
	\geq 
		\Pr\left(E_K\right) - \frac{c_0 \stepsize_K}{\epsilon} 
	= 1 - \Pr\left( s \leq K\right) - \frac{c_0 \stepsize_K}{\epsilon}.
\end{align} 

Finally, we will upper bound $\Pr\left( s \leq K\right)$.  For $k=0, 1, \ldots$, define the variables $n_k = \min(k,s)$, we then have
\begin{align}\label{eq:thrm:cvx:minibatch:finite:proof:Pne}
	\Pr\left( s \leq K\right) 
	= \Pr\left( \max\{e_0, \ldots, e_K\} > M\right) 
	&= \Pr\left( e_{n_s}^2 \bindic{s \leq K}  > M^2 \right) 
	\nonumber\\
	&\mathop \leq \limits^{\mathrm{(a)}} 	 \E\left[e_{n_s}^2\bindic{s \leq K}  \right] / M^2
	\nonumber\\
	&= \E\left[e_{n_K}^2\bindic{s \leq K}  \right] / M^2	
	\nonumber\\
	&\leq \E\left[e_{n_K}^2 \right] / M^2,
\end{align}
where $\mathrm{(a)}$ follows from Markov's inequality.
To bound $ \E\left[e_{n_K}^2\right]$, let $\xi_k := e_{n_{k}}^2$, we have
\begin{align*}
	\E \left[\xi_{k+1} \big | \mc{F}_k\right]
	&\leq
		\E \left[\xi_{k+1} \bindic{s \leq k}  \big | \mc{F}_k\right]
		+
		\E \left[\xi_{k+1}  \bindic{s > k} \big | \mc{F}_k\right]
	\nonumber\\
	&=
		\xi_{k} \bindic{s \leq k}  
		+
		 \E \left[e_{k+1}^2 \bindic{E_k}  \big | \mc{F}_k\right]
	\nonumber\\
	&\leq 
		\xi_{k} \bindic{s \leq k}  
		+
	    e_k^2 \bindic{E_k}\left(1 - \mu \varrho\stepsize_k\right)  + 
	     \left({ \sigma^2  }/{\mu } +  \gamma^2\right)\stepsize_k^2
	  \nonumber\\
	  &=
	  	\xi_{k} \bindic{s \leq k}  
	  		+
	  	    \xi_k \bindic{s> k} \left(1 - \mu \varrho\stepsize_k\right) + 
	  	     \left({ \sigma^2  }/{\mu } +  \gamma^2\right)\stepsize_k^2
	 \nonumber\\
	 &\leq
		 \xi_k +  \left({ \sigma^2  }/{\mu } +  \gamma^2\right)\stepsize_k^2,
\end{align*}
where the second inequality follows from \eqref{eq:thrm:cvx:minibatch:finite:proof:e:1} and the fact that $m_k=1/\stepsize_k$.
We successively deduce that
\begin{align*}
	\E \left[e_{n_K}\right] = \E \left[\xi_{K} \right]
	\leq 
		\xi_0 + \left({ \sigma^2  }/{\mu } +  \gamma^2\right) \sum_{k=0}^{K-1} \stepsize_k^2.
\end{align*}
Plugging the preceding result into \eqref{eq:thrm:cvx:minibatch:finite:proof:Pne} yields
\begin{align}\label{eq:thrm:cvx:minibatch:finite:proof:Pne:2}
	\Pr\left( s \leq K\right) 
	\leq
		\frac{\xi_0 +  \left({ \sigma^2  }/{\mu } +  \gamma^2\right) \sum_{k=0}^{K-1} \stepsize_k^2 }{M^2} .
\end{align}
Finally, by selecting $M = \sqrt{\xi_0/\delta}$ and combining \eqref{eq:thrm:cvx:minibatch:finite:proof:Pe:2} and \eqref{eq:thrm:cvx:minibatch:finite:proof:Pne:2} give the first claim in the theorem.

For the second claim, we first notice that the stepsizes $\stepsize_k$ are now constant and that $\mu\varrho \stepsize \in (0,1)$ since $\stepsize_0\leq 1/(\mu\varrho)$.  It is thus easy to verify from \eqref{eq:thrm:cvx:minibatch:finite:proof:e:2} that
\begin{align*}
	\E\left[ e_{K}^2 \bindic{E_{K}}\right] 
	\leq
		\left(1-\mu\varrho\stepsize\right)^K e_0^2 
		+
		\frac{\left( \sigma^2 /  \mu +  \gamma^2\right)\stepsize}{\mu \varrho}.
%	\leq 
%		 \frac{2\left( \sigma^2 /  \mu +  \gamma^2\right)\stepsize}{\mu \varrho} 
\end{align*}
Let $\eta={\left( \sigma^2 /  \mu +  \gamma^2\right)}/{(\mu \varrho)}$, we  have $\left(1-\mu\varrho\stepsize\right)^K e_0^2  \leq \eta \stepsize$ whenever $K\geq \frac{1}{\mu\varrho\stepsize}\log\left(\frac{e_0^2}{\eta\stepsize} \right)$. 
Now since $\stepsize= \stepsize_0 K^{-\tau}$ with $\tau \in (1/2,1)$, taking $K$ be such that
\begin{align*}
	K \geq \left(\frac{\log\left( \frac{e_0^2 K^{\tau}}{\eta \stepsize_0}\right)}{\mu\varrho\stepsize_0}\right)^{\frac{1}{1-\tau}},
\end{align*}
suffices to guarantee $K\geq \frac{1}{\mu\varrho\stepsize}\log\left(\frac{e_0^2}{\eta\stepsize} \right)$. In that case, we have $
	\E\left[ e_{K}^2 \bindic{E_{K}}\right] 
	\leq
	2 \eta \stepsize
$. Consequently, replacing the constant $c_0$ in the first claim by $2\eta$ and setting $\epsilon=\frac{2\eta\stepsize_0}{\delta K^{\tau}}$ completes the proof.

\section{Proof of Lemma~\ref{lem:xdiff:grad:bound}}\label{appendix:proof:lem:xdiff:grad:bound}
{Let $\Delta_i := \ltwo{f'(x_i) - g_i }^2$ and $x^\star = \gproj{\Xopt}{x_k}$,  it follows that from equation~\eqref{eq:proposition:stability:proof:xdiff:2} that
\begin{align*} 
	\ltwo{x_{k}-x\opt} ^2
	\leq 
		 \dist(x_0, \Xopt) ^2
		+ \sum_{i=0}^{k-1} \stepsize_i\left(\Delta_i/(2\mu) + \gamma^2\right).
\end{align*}
Let $q\geq 2$, we deduce that 
\begin{align*} 
	\E\left[\ltwo{x_{k}-x\opt} ^ q \right]
	&\leq 
		\E\left[\left(
			\dist(x_0, \Xopt) ^2 
			+ 
			\sum_{i=0}^{k-1} \stepsize_i\left(\Delta_i/(2\mu) + \gamma^2\right)\right)^{\frac{q}{2}}
		\right]
	\nonumber\\
	&\leq
		2^{\frac{q}{2}} \left(\dist(x_0, \Xopt) ^2 + \gamma^2\sum_{i=0}^{k} \stepsize_i\right)^{\frac{q}{2}}
		+  
	 		\mu^{-\frac{q}{2}}  \E\left[\left(\sum_{i=0}^{k-1} \stepsize_i\Delta_i\right)^{\frac{q}{2}} 	\right]	 	
	 \nonumber\\
	 &\leq
	 		2^{q} \dist(x_0, \Xopt) ^q + (2\gamma)^{q} \left(\sum_{i=0}^{k-1} \stepsize_i\right)^{\frac{q}{2}}
%	 \nonumber\\
%	 &\hspace{0.45cm}
	 		+
	 		\mu^{-\frac{q}{2}} 
	 		\E\left[\left(\sum_{i=0}^{k-1} \frac{\stepsize_i }{\sum_{i=0}^{k-1} \stepsize_i}\Delta_i\right)^{\frac{q}{2}} 	\right] \left(\sum_{i=0}^{k-1} \stepsize_i\right)^{\frac{q}{2}} ,
\end{align*}
where we applied the  inequality $(a+b)^\epsilon \leq 2^\epsilon( a^\epsilon +  b^\epsilon)$ for any $a,b, \epsilon \geq 0$ twice, as well as multiplied and divided the same quantity in the last step.
Using Jensen's inequality and Assumption~\ref{assumption:polynomial:growth}, the right most term in the preceding inequality can be upper bounded by
\begin{align*}
	 	\sum_{i=0}^{k-1} \frac{\stepsize_i }{\sum_{i=0}^{k-1} \stepsize_i}\E\left[\Delta_i\right]^{\frac{q}{2}} \left(\sum_{i=0}^{k-1} \stepsize_i\right)^{\frac{q}{2}} 
	 	\leq 
	 	\sigma^{\frac{q}{4}+1} \left(\sum_{i=0}^{k-1} \stepsize_i\right)^{\frac{q}{2}}.
\end{align*}
Since $\stepsize_i = \stepsize_0(i+1)^{-\tau}$ with $\tau \in (0,1)$, we have for $k\geq 1$:
\begin{align*}
	\sum_{i=0}^{k-1} \stepsize_i = \sum_{j=1}^{k} \stepsize_0 j^{-\tau} 
	\leq \stepsize_0\left(1 + \int_{1}^{k} t^{-\tau}\mathrm{d}t \right)
	\leq \stepsize_0\left(1+\frac{ k^{1-\tau}}{1-\tau}\right)
	\leq \frac{2\stepsize_0}{1-\tau}k^{1-\tau}.
\end{align*}
Collecting the terms, we obtain
\begin{align*} 
	\E\left[\ltwo{x_{k}-x\opt} ^ q \right]
	\leq 
		2^{q}\dist(x_0, \Xopt) ^q
		+ 
		\left(  (2\gamma)^q + \mu^{-\frac{q}{2}} 	\sigma^{\frac{q}{4}+1} \right) \left(\frac{2\stepsize_0}{1-\tau}k^{1-\tau}\right)^{\frac{q}{2}}.
\end{align*}			
Therefore, setting $q=4(p-1)$ gives the second bound in the lemma. 
Finally, setting $q=2(p-1)$ and using Assumption~\ref{assumption:polynomial:growth} yields the first bound in the lemma since
\begin{align*}
	\E\left[\ltwo{f'(x_k, S)}^2\right] 
	= 
		\E\left[\E\left[\ltwo{f'(x_k, S)}^2 \big| \mc{F}_k\right]\right]
		 \leq  L_0 + L_1\E\left[\dist\left(x_k, \Xopt\right)^{2(p-1)}\right].
\end{align*}}

\section{Proof of Theorem~\ref{thrm:cvx:finite:tvclip}}
We rely on the following useful lemma \cite[Lemma~4]{Chu54}, see also \cite[Lemma~5, Page~46]{Pol87}.
\begin{lemma}[Chung's lemma]\label{lem:chung}
Let $V_k$ be a sequence of nonnegative random variables satisfying:
\begin{align*}
	V_{k+1}  \leq \left(1-\frac{a}{(k+m)^\tau}\right) V_k + \frac{b}{(k+m)^{t}}, \quad 0 < \tau < 1, \quad t > \tau,  \quad m \in \N_{+},
\end{align*}
then 
\begin{align*}
		V_k \leq 	\frac{b}{a}\frac{1}{k^{t-\tau}} + o\left(\frac{1}{k^{t-\tau}}\right).
\end{align*}
\end{lemma}
\medskip

\noindent We begin by letting $\varrho_k = \min\left\{1, \gamma_k/\ltwo{g_k}\right\} $, $x\opt = \gproj{\Xopt}{x_k}$, and $f'(x_k) = \E[g_k | \mc{F}_k]$. It holds  that
\begin{align*}
	 \ltwo{x_{k+1}-x\opt} ^2
%	&= 
%		 \ltwo{x_k-x\opt} ^2
%		- 2\stepsize_k \InP{d_k}{x_k-x\opt}
%		+ \stepsize_k^2 \ltwo{d_k}^2
%	\nonumber\\
	&=
			 \ltwo{x_k-x\opt} ^2
			- 2\stepsize_k  \InP{f'(x_k)}{x_k-x\opt}
			+ 2 \stepsize_k \InP{ f'(x_k) - \varrho_k g_k}{x_k-x\opt}
			+ \stepsize_k^2 \ltwo{d_k}^2.
\end{align*}
Again by the convexity of $f$ and Assumption~\ref{assumption:quadratic:growth}, we have $
 \InP{f'(x_k)}{x_k-x\opt} \geq \mu \ltwo{x_k-x\opt}^2.$ Note also that $\ltwo{d_k} = \ltwo{\varrho_k g_k} \leq \ltwo{g_k}$ since $\varrho_k\leq 1$. Therefore, plugging these results into the preceding equation and taking  the conditional expectation gives
\begin{align}\label{eq:thrm:cvx:rate:proof:xdiff:1}
	\E \left[\ltwo{x_{k+1}-x\opt} ^2  \big | \mc{F}_k\right] 
	&\leq
			\left(1- 2\mu\stepsize_k\right) \ltwo{x_k-x\opt} ^2 
			+ 2 \stepsize_k \InP{ \Delta_k }{x_k-x\opt}
			+ \stepsize_k^2 \E\left[\ltwo{g_k}^2 \big | \mc{F}_k\right]
	\nonumber\\
	&\leq
			\left(1- \mu\stepsize_k\right) \ltwo{x_k-x\opt} ^2 
			+ \frac{\stepsize_k}{\mu} \ltwo{\Delta_k}^2
			+ \stepsize_k^2 \E\left[\ltwo{g_k}^2 \big | \mc{F}_k\right],
\end{align}
where $\Delta_k=f'(x_k) - \E\left[\varrho_k g_k \big | \mc{F}_k\right]$, and the second inequality follows from Young's inequality \eqref{eq:young} with $\epsilon=\mu$. The term $\ltwo{\Delta_k}$ can be bounded as
\begin{align*}
\ltwo{ \Delta_k }
	= \ltwo{ \E\left[ g_k - \varrho_k g_k \big | \mc{F}_k \right]}
	&\leq  \ltwo{ \E\left[ g_k \bindic{\ltwo{g_k}\geq \gamma_k} \big | \mc{F}_k \right]}
	\nonumber\\
	&\leq
	\E \left[ \ltwo{g_k \bindic{\ltwo{g_k}\geq \gamma_k} } \big | \mc{F}_k  \right]
	\nonumber\\
	&\leq \E\left[ \ltwo{  g_k}^2 \gamma_k^{-1} \big | \mc{F}_k \right],
\end{align*}
where the first inequality follows from Jensen's inequality. Now, since by Assumption~\ref{assumption:polynomial:growth}, $\E\left[ \ltwo{  g_k}^2 \big | \mc{F}_k \right] \leq L_0 + L_1 \dist\left(x_k, \Xopt\right)^{2(p-1)}$, we deduce that
\begin{align}\label{eq:thrm:cvx:rate:proof:noise:1}
	\ltwo{\Delta_k}^2 
	\leq 
		\gamma_k^{-2} \left(L_0 + L_1\dist\left(x_k, \Xopt\right)^{2(p-1)}\right)^2
%	\nonumber\\
	\leq 
		2\gamma_k^{-2} \left(L_0^2 + L_1^2\dist\left(x_k, \Xopt\right)^{4(p-1)}\right),
\end{align}
where we used the inequality $(a+b)^2 \leq 2a^2 + 2b^2$ in the last step.
Plugging \eqref{eq:thrm:cvx:rate:proof:noise:1} into \eqref{eq:thrm:cvx:rate:proof:xdiff:1} and taking the full expectation of the resulting expression yields
\begin{align*}
	\E \left[\ltwo{x_{k+1}-x\opt} ^2\right] 
	&\leq
			\left(1- \mu\stepsize_k\right) 	\E \left[\ltwo{x_{k}-x\opt} ^2\right] 
	\nonumber\\
	&\hspace{0.45cm}
			+  \frac{2\stepsize_k \gamma_k^{-2}}{\mu} \left(L_0^2 + L_1^2  \E\left[\dist\left(x_k, \Xopt\right)^{4(p-1)} \right]\right)
			+ \stepsize_k^2 \E\left[\ltwo{g_k}^2\right].
\end{align*}
Now, we can invoke Lemma~\ref{lem:xdiff:grad:bound} to obtain
\begin{align}\label{eq:thrm:cvx:rate:proof:xdiff:2}
	\E \left[\ltwo{x_{k+1}-x\opt} ^2\right] 
	\leq
			\left(1- \mu\stepsize_k\right) 	\E \left[\ltwo{x_{k}-x\opt} ^2\right] 
			+  C_0 \stepsize_k \gamma_k^{-2} k^{2(p-1)(1-\tau)}
			+ C_1\stepsize_k^2  k^{(p-1)(1-\tau)},
\end{align}
where $C_0 = \frac{2}{\mu}(L_0^2 + L_1^2 (D_0+D_1))$ and $C_1 = G_0+G_1$. Finally, since $\stepsize_k= \frac{\stepsize_0}{(k+1)^\tau}$ and $\gamma_k = \gamma/\sqrt{\stepsize_k}$, we get
\begin{align*}
	\E \left[\ltwo{x_{k+1}-x\opt} ^2 \right] 
	&\leq
			\left(1- \frac{\mu\stepsize_0}{(k+1)^\tau}\right) \E \left[\ltwo{x_k-x\opt} ^2 \right]
			+ \frac{C}{(k+1)^{2 \left(1- p(1-\tau)\right)}},
\end{align*}
where $C = C_0\gamma^2+ C_1$. Finally, plugging $\tau = 1- \epsilon$ with $\epsilon \in (0,1)$ into the preceding inequality and invoking Lemma~\ref{lem:chung} completes the proof.

\section{Proof of Lemma~\ref{lem:xdiff:non-smooth}}

We first recall  the update formula  of the search direction $d_k$:
\begin{align*}%\label{eq:alg:dk}
	d_k &= \clip{(1-\mmt_{k-1}) d_{k-1} + \mmt_{k-1} g_k}
%	\nonumber\\
%	&= 
	=
	\argmin_{z\in \mathrm{B}(0, \gamma)} \left\{ \frac{1}{2} \ltwo{z - (1-\mmt_{k-1}) d_{k-1} - \mmt_{k-1} g_k}^2\right\},
\end{align*}
whose optimality condition  implies that
\begin{align*}
	\InP{(1-\mmt_{k-1}) (d_{k-1}-d_k) + \mmt_{k-1} ( g_k -d_k)}{ z - d_k} \leq 0 \quad \forall z\in \mathrm{B}(0, \gamma).
\end{align*}
We deduce (by taking $z=0$) that
\begin{align}\label{eq:lem:dk:proof:dksq}
	\ltwo{d_k}^2 
	&\leq 
		\mmt_{k-1}\InP{g_k}{d_k}  + (1-\mmt_{k-1}) \InP{d_{k-1}}{d_k}
	\nonumber\\
	&=			
		\mmt_{k-1}\InP{g_k}{d_k} 
		+ 
		(1-\mmt_{k-1}) \left(\frac{1}{2}\ltwo{d_{k-1}}^2 + \frac{1}{2}\ltwo{d_k}^2 - \frac{1}{2}\ltwo{d_k-d_{k-1}}^2\right).
\end{align}
We have
\begin{align}\label{eq:lem:dk:proof:g:ddiff}
	 &\mmt_{k-1}\InP{g_k}{d_k} - \frac{1-\mmt_{k-1}}{2}\ltwo{d_k-d_{k-1}}^2 
	\nonumber\\
	&\hspace{2cm}		 
	 = 
	 	 \mmt_{k-1}\InP{g_k}{d_{k-1}} 
	 	+ \mmt_{k-1}\InP{g_k}{  d_k - d_{k-1}  } 
	 	- \frac{1-\mmt_{k-1}}{2}\ltwo{d_k-d_{k-1}}^2
	 \nonumber\\
	 &\hspace{2cm}
	 \leq
	 	 \mmt_{k-1}\InP{g_k}{d_{k-1}} 
	 	+
	 	\frac{\mmt_{k-1}^2}{2(1-\mmt_{k-1})}\ltwo{g_k}^2,
\end{align}
where the last step follows from Young's inequality \eqref{eq:young} with $\epsilon=(1-\mmt_{k-1})/\mmt_{k-1}$.
Plugging \eqref{eq:lem:dk:proof:g:ddiff} into \eqref{eq:lem:dk:proof:dksq} and rearranging yields
\begin{align}\label{eq:lem:dk:proof:key:1}
	\frac{1}{2}\ltwo{d_k}^2 
	\leq			
		\frac{1}{2}\ltwo{d_{k-1}}^2
		+ \mmt_{k-1}\InP{g_k}{d_{k-1}} 
		- \frac{\mmt_{k-1}}{2}\left(\ltwo{d_{k-1}}^2 + \ltwo{d_k}^2\right) 
		+ 	\frac{\mmt_{k-1}^2}{2(1-\mmt_{k-1})}\ltwo{g_k}^2.
\end{align}
By the weak convexity of $f$ and Assumption~(A1), it holds that
\begin{align}\label{eq:lem:dk:proof:key:2}
	\mmt_{k-1} \,\E[\InP{g_k}{d_{k-1}} | \mc{F}_{k}] 
	&= \nu \InP{\E[g_k | \mc{F}_{k}] }{x_{k-1}-x_k} 
	\nonumber\\
	&\leq \nu\left(f(x_{k-1}) -f(x_k) + \frac{\rho}{2} \ltwo{x_{k-1}-x_k}^2\right).
\end{align}
We can now take the conditional expectation in \eqref{eq:lem:dk:proof:key:1}, combine the result with \eqref{eq:lem:dk:proof:key:2}, and rearrange terms to obtain
\begin{align*}
	\nu f(x_k)	+ \E\left[\frac{1}{2} \ltwo{d_k}^2  \bigg| \mc{F}_{k}\right]
	&\leq
		\nu f(x_{k-1} ) + \frac{1-\mmt_{k-1}}{2} \ltwo{d_{k-1}}^2 
%	\nonumber\\
%	&\hspace{0.45cm}
		-\frac{\mmt_{k-1}}{2} \E[ \ltwo{d_k}^2  \big| \mc{F}_{k}] 
	\nonumber\\
	&\hspace{0.45cm}
		+ \frac{\mmt_{k-1}^2}{2(1-\mmt_{k-1})} \E[\ltwo{g_k}^2  \big| \mc{F}_{k}]  + \frac{\rho \nu}{2}\ltwo{x_k-x_{k-1}}^2.
\end{align*}
Since $\mmt_k \leq \mmt_{k-1}$, it follows that
\begin{align*}
	\nu f(x_k)	+ \E\left[\frac{1-\mmt_k}{2}  \ltwo{d_k}^2  \bigg | \mc{F}_{k}\right]
	&\leq
		\nu f(x_{k-1} ) + \frac{1-\mmt_{k-1}}{2} \ltwo{d_{k-1}}^2 
%	\nonumber\\
%	&\hspace{0.45cm}
		-\mmt_{k} \E[ \ltwo{d_k}^2  \big| \mc{F}_{k}]
	\nonumber\\
	&\hspace{0.45cm}
		+ \frac{\mmt_{k-1}^2}{2(1-\mmt_{k-1})} \E[\ltwo{g_k}^2  \big| \mc{F}_{k}]  + \frac{\rho \nu}{2}\ltwo{x_k-x_{k-1}}^2.
\end{align*}
By Assumption~\ref{assumption:wcvx:gradient:boundedness}, we have $\E[\ltwo{g_k}^2  \big| \mc{F}_{k}]\leq L^2$. Note also that $\ltwo{x_{k}-x_{k-1}}^2 	=	\ltwo{\stepsize_{k-1} d_{k-1}}^2 \leq \stepsize_{k-1}^2 \thres^2$. 
Therefore, multiplying both sides of the preceding inequality by $1/\nu = \stepsize_{k-1}/\mmt_{k-1}=\stepsize_k/\mmt_k$ completes the proof.

\section{Proof of Lemma~\ref{lem:wcvx:Wfunc}}

We start by showing that the gradient of $f_\lambda$ is bounded. 
Recall that
\begin{align*}
	\prox{\lambda f}{x} = \argmin_{y\in\R^n} \left\{f(y)+ \frac{1}{2\lambda} \ltwo{x-y}^2\right\},
\end{align*}
which implies for any $x\in \R^n$ that
\begin{align*}
	\frac{1}{2\lambda}\ltwo{x - \prox{\lambda f}{x} }^2 
	\leq 
		f(x) - f(\prox{\lambda f}{x} ).
\end{align*}
Since the function $f$, satisfying Assumption~\ref{assumption:wcvx:gradient:boundedness}, is necessarily Lipschitz continuous with constant $L$ \cite{DD19}, it holds that  $f(x) - f(\prox{\lambda f}{x} ) \leq L \ltwo{x- \prox{\lambda f}{x} }$. Thus, combining these results and the definition of  $\grad{f_{\lambda}}$ yields
$\ltwo{\grad{f_{\lambda}(x_k)}} \leq 2L \leq \gamma$. We deduce that  $\clip{\grad{f_{\lambda}(x_k)}} = \grad{f_{\lambda}(x_k)}$.\medskip

 For simplicity, let $z_k=(1-\mmt_{k-1}) d_{k-1} + \mmt_{k-1} g_k$, and hence $d_k=\clip{z_k}$. We can now proceed to bound the quantity $\half\ltwo{d_k - \grad{f_{\lambda}(x_k)}}^2$ as follows:
\begin{align}\label{eq:wcvx:lem:2:proof:direction:1}
	\half\ltwo{d_k - \grad{f_{\lambda}(x_k)}}^2
	&=
			\half\ltwo{\clip{z_k} - \clip{\grad{f_{\lambda}(x_k)}} }^2
	\leq
		\half\ltwo{z_k - \grad{f_{\lambda}(x_k)}}^2
	\nonumber\\
	&= 
		\half\ltwo{z_k }^2
		- 
		\InP{z_k}{ \grad{f_{\lambda}(x_k)}}
		+				
		\half\ltwo{\grad{f_{\lambda}(x_k)}}^2,
\end{align}
where the inequality follows from the non-expansiveness of the clipping operator.
Next, we develop the middle term in \eqref{eq:wcvx:lem:2:proof:direction:1} as
\begin{align}\label{eq:wcvx:lem:2:proof:inp:1}
	-\InP{z_k}{ \grad{f_{\lambda}(x_k)}}
%	&= 
%		-\InP{(1-\mmt_{k-1}) d_{k-1} + \mmt_{k-1} g_k}{\grad{f_{\lambda}(x_k)}}
%	\nonumber\\
	&=
		-\mmt_{k-1} \InP{g_k}{\grad{f_{\lambda}(x_k)}}
		-(1-\mmt_{k-1})\InP{d_{k-1}}{\grad{f_{\lambda}(x_{k-1})}}
	\nonumber\\		
	&\hspace{0.45cm}	
		+ (1-\mmt_{k-1}) \InP{d_{k-1}}{\grad{f_{\lambda}(x_{k-1})} - \grad{f_{\lambda}(x_k)}}
	\nonumber\\
	&=
		-\mmt_{k-1}\InP{g_k}{\grad{f_{\lambda}(x_k)}} 
		+\mmt_{k-1} \InP{d_{k-1}}{\grad{f_{\lambda}(x_{k-1})}}
	\nonumber\\	
	&\hspace{0.45cm}	
	+	
		\half\ltwo{d_{k-1} - \grad{f_{\lambda}(x_{k-1})}}^2 
		-\half\ltwo{d_{k-1}}^2 
		-\half\ltwo{\grad{f_{\lambda}(x_{k-1})}}^2
	\nonumber\\		
	&\hspace{0.45cm}	
		+ (1-\mmt_{k-1}) \InP{d_{k-1}}{\grad{f_{\lambda}(x_{k-1})} - \grad{f_{\lambda}(x_k)}}.
\end{align}
Plugging \eqref{eq:wcvx:lem:2:proof:inp:1} into \eqref{eq:wcvx:lem:2:proof:direction:1} and rearranging terms yields
\begin{align}\label{eq:wcvx:lem:2:proof:direction:2}
	\half\ltwo{d_k - \grad{f_{\lambda}(x_k)}}^2 - 	\half\ltwo{\grad{f_{\lambda}(x_k)}}^2	
	&\leq
		\half\ltwo{d_{k-1} - \grad{f_{\lambda}(x_{k-1})}}^2 - \half\ltwo{\grad{f_{\lambda}(x_{k-1})}}^2
	\nonumber\\		
	&\hspace{0.45cm}
		-\mmt_{k-1} \InP{g_k}{\grad{f_{\lambda}(x_k)}} 
		+\mmt_{k-1} \InP{d_{k-1}}{\grad{f_{\lambda}(x_{k-1})}}
	\nonumber\\		
	&\hspace{0.45cm}		
		+
			\half\ltwo{z_k }^2 -\half\ltwo{d_{k-1}}^2 
	\nonumber\\		
	&\hspace{0.45cm}	
	+ (1-\mmt_{k-1}) \InP{d_{k-1}}{\grad{f_{\lambda}(x_{k-1})} - \grad{f_{\lambda}(x_k)}}.
\end{align}
Using the definition of $z_k$, it is easy to verify that
\begin{align*}
	\half \ltwo{z_k}^2 
%	=
%		\half 
%		\ltwo{(1-\mmt_{k-1})d_{k-1} + \mmt_{k-1} g_k}^2 
	-
		\half \ltwo{d_{k-1}}^2 
	&=
		\frac{\mmt_{k-1}^2}{2} \ltwo{d_{k-1} - g_k}^2 
		+
		\mmt_{k-1} \InP{g_k-d_{k-1}}{d_{k-1} }
	\nonumber\\
	&\leq 
		\mmt_{k-1}^2 \ltwo{d_{k-1}}^2 
		+
		\mmt_{k-1}^2 \ltwo{g_k}^2 
		+
		\mmt_{k-1}\InP{g_k}{d_{k-1} }
		-\mmt_{k-1}\ltwo{d_{k-1}}^2
	\nonumber\\
	&\leq
		\mmt_{k-1}^2 \ltwo{g_k}^2 
		+
		\nu \InP{g_k}{x_{k-1} - x_k },	
\end{align*}
where we used $\mmt_{k-1}^2\leq \mmt_{k-1}$ in the last step.
Consequently, using the same arguments leading to eq.~\eqref{eq:lem:dk:proof:key:2} in the proof of Lemma~\ref{lem:xdiff:non-smooth} and Assumption~\ref{assumption:wcvx:gradient:boundedness}, we get 
\begin{align}\label{eq:wcvx:lem:2:proof:direction:zdiff}
	\E\left[\half \ltwo{z_k}^2  \big | \mc{F}_k \right]
	-
	\half \ltwo{d_{k-1}}^2 
	&\leq
		\mmt_{k-1}^2\E\left[  \ltwo{g_k}^2 \big | \mc{F}_k \right]
		+
			\nu \InP{\E\left[g_k \big | \mc{F}_k \right]}{x_{k-1} - x_k } 
	\nonumber\\
	&\leq 
		\mmt_{k-1}^2 L^2
		+
		 \nu\left(f(x_{k-1}) -f(x_k) + \frac{\rho}{2} \ltwo{x_{k-1}-x_k}^2\right)
	\nonumber\\
	&\leq 
		\mmt_{k-1}^2 L^2
		+
		 \nu\left(f(x_{k-1}) -f(x_k) + {\stepsize_{k-1}^2\rho\thres^2}/{2} \right),
\end{align}
where the last inequality holds since $\ltwo{x_{k}-x_{k-1}}^2 	=	\ltwo{\stepsize_{k-1} d_{k-1}}^2 \leq \stepsize_{k-1}^2 \thres^2$. 
Since $f_{\lambda}$ is $(1/\lambda)$-smooth, the right most term in \eqref{eq:wcvx:lem:2:proof:inp:1} can be bounded as
\begin{align}\label{eq:wcvx:lem:2:proof:direction:graddiff}
	(1-\mmt_{k-1}) \InP{d_{k-1}}{\grad{f_{\lambda}(x_{k-1})} - \grad{f_{\lambda}(x_k)}}
	\leq 
		\frac{\stepsize_{k-1} (1-\mmt_{k-1}) }{\lambda } \ltwo{d_{k-1}}^2 
	\leq 
			\frac{\stepsize_{k-1}}{\lambda } \ltwo{d_{k-1}}^2 .
\end{align}
Therefore, by taking the conditional expectation in \eqref{eq:wcvx:lem:2:proof:direction:2} and combining the result with \eqref{eq:wcvx:lem:2:proof:direction:zdiff} and \eqref{eq:wcvx:lem:2:proof:direction:graddiff}, we arrive at
\begin{align}\label{eq:wcvx:lem:2:proof:direction:3}
	&\E\left[\half\ltwo{d_k - \grad{f_{\lambda}(x_k)}}^2 \big | \mc{F}_k\right] - 	\half\ltwo{\grad{f_{\lambda}(x_k)}}^2	+ 	\nu f(x_k)
	\nonumber\\
	&\hspace{2.0cm}\leq
		\half\ltwo{d_{k-1} - \grad{f_{\lambda}(x_{k-1})}}^2 - \half\ltwo{\grad{f_{\lambda}(x_{k-1})}}^2
		+ \nu f(x_{k-1})
	\nonumber\\		
	&\hspace{2.45cm}
		-\mmt_{k-1} 	\E\left[ \InP{g_k}{\grad{f_{\lambda}(x_k)}} \big | \mc{F}_k\right] 
		+\mmt_{k-1} \InP{d_{k-1}}{\grad{f_{\lambda}(x_{k-1})}}
	\nonumber\\		
	&\hspace{2.45cm}		
		+\frac{\stepsize_{k-1}}{\lambda } \ltwo{d_{k-1}}^2
		+\mmt_{k-1}^2 L^2 + {\nu \stepsize_{k-1}^2\rho\thres^2}/{2}.
\end{align}
Finally, multiplying both sides of \eqref{eq:wcvx:lem:2:proof:direction:3} by $1/\nu = \stepsize_{k-1}/\mmt_{k-1}$ completes the proof.

\section{Proof of Lemma~\ref{lem:Vfunc:non-smooth}}

We begin with the smoothness of $f_\lambda$ which ensures (see, Lemma~\ref{lem:moreau:env}) that
\begin{align}\label{eq:wcvx:lem:2:proof:envelope:1}
	f_\lambda(x_k) 
	&\leq 
		f_\lambda(x_{k-1}) 
		+ \InP{\grad{f_{\lambda}(x_{k-1})}}{x_k-x_{k-1}} 
		+ \frac{1}{2\lambda}\ltwo{x_k-x_{k-1}}^2
	\nonumber\\		
	&=
		f_\lambda(x_{k-1}) 
		-\stepsize_{k-1} \InP{\grad{f_{\lambda}(x_{k-1})}}{d_{k-1}}
		+ \frac{1}{2\lambda}\ltwo{x_k-x_{k-1}}^2.		
\end{align}
Since $\ltwo{x_{k}-x_{k-1}}^2 	=	\ltwo{\stepsize_{k-1} d_{k-1}}^2 \leq \stepsize_{k-1}^2 \thres^2$ and $\rho\leq 1/(2\lambda)$,  
summing \eqref{eq:wcvx:lem:2:proof:envelope:1} and \eqref{eq:lem:Wk} yields
\begin{align}\label{eq:wcvx:lem:2:proof:envelope:2}
	f_\lambda(x_k) + \E\left[W_k\big| \mc{F}_k\right]
	&\leq
		 f_\lambda(x_{k-1}) + W_{k-1}	
		-\stepsize_{k-1} 	\E\left[ \InP{g_k}{\grad{f_{\lambda}(x_k)}} \big | \mc{F}_k\right] 
	\nonumber\\		
	&\hspace{0.45cm}	+\frac{\stepsize_{k-1}}{\lambda\nu } \ltwo{d_{k-1}}^2
		+ \stepsize_{k-1}^2 \left(\nu L^2 + \thres^2/\lambda \right).		 
\end{align}
Next, we bound the term $\E\left[ \InP{g_k}{\grad{f_{\lambda}(x_k)}} \big | \mc{F}_k\right] $. 
Consider the proxial point $\hat{x}_{k+1}$ defined as
\begin{align*}
	\hat{x}_{k+1}:=\prox{\lambda f}{x_k} = \argmin_{y\in\R^n} \left\{f(y) +  \frac{1}{2\lambda}\ltwo{x_k-y}^2\right\},
\end{align*}
it then follows from  the weak convexity of $f$ and the definition of $\grad{f_{\lambda}(x_k)}$ that
\begin{align*}
	-\E [\InP{g_k}{\grad{f_{\lambda}(x_k)}} | \mc{F}_{k} ]
	&=
		\lambda^{-1} \InP{\E[g_k | \mc{F}_k]}{\hat{x}_{k+1} - x_k } 
	\nonumber\\			
	&\leq
		\lambda^{-1}
		\left[ f(\hat{x}_{k+1}) - f(x_k) + \frac{\rho}{2}\ltwo{\hat{x}_{k+1}-x_k}^2 \right].
\end{align*}
Since the function $x \mapsto f(x) + \frac{1}{2\lambda}\ltwo{x-x_k}^2$ \, is  $(\lambda^{-1}-\rho)$-strongly convex with $\hat{x}_{k+1}$ being its minimizer, it follows from \cite[Theorem~5.25]{Bec17} that
\begin{align}\label{eq:strong:cvx}
	f(x_k) + \frac{1}{2\lambda} \ltwo{x_k-x_k}^2
	-
	\left(f(\hat{x}_{k+1}) + \frac{1}{2\lambda} \ltwo{\hat{x}_{k+1}-x_k}^2\right)
	\geq
		\frac{\lambda^{-1}-\rho}{2}\ltwo{\hat{x}_{k+1}-x_k}^2.
\end{align}		
 We thus have
 \begin{align}\label{eq:proof:thrm:1:inp:1:bounded}
  f(\hat{x}_{k+1}) - f(x_k) + \frac{\rho}{2}\ltwo{\hat{x}_{k+1}-x_k}^2
 	&=
 		f(\hat{x}_{k+1}) + \frac{1}{2\lambda} \ltwo{\hat{x}_{k+1}-x_k}^2
 		-
 		f(x_k)
 	\nonumber\\
 	&\hspace{0.45cm} 	
 		-
 		\frac{\lambda^{-1}-\rho}{2} \ltwo{\hat{x}_{k+1}-x_k}^2
 	\nonumber\\
 	&\hspace{0.0cm}
 	\mathop \leq \limits^{\mathrm{(a)}} 	
 		-
 		(\lambda^{-1}-\rho) \ltwo{\hat{x}_{k+1}-x_k}^2
 	\nonumber\\
 	&\hspace{.0cm}
 	\mathop = \limits^{\mathrm{(b)}}
 		-\lambda^2(\lambda^{-1}-\rho) \ltwo{\grad{f_{\lambda}(x_k)}}^2,
 \end{align}
where $\mathrm{(a)}$ is due to \eqref{eq:strong:cvx} and  $\mathrm{(b)}$ follows from the definition of $\grad{f_{\lambda}(x_{k})}$.
Plugging \eqref{eq:proof:thrm:1:inp:1:bounded} into \eqref{eq:wcvx:lem:2:proof:envelope:2} yields
\begin{align}\label{eq:wcvx:lem:2:proof:envelope:3}
	f_\lambda(x_k) + \E\left[W_k\big| \mc{F}_k\right]
	&\leq
		 f_\lambda(x_{k-1}) + W_{k-1}	+ \frac{\stepsize_{k-1}}{\lambda\nu } \ltwo{d_{k-1}}^2
	\nonumber\\		
	&\hspace{0.45cm}	
		-\stepsize_{k-1} (1-\rho\lambda) \ltwo{\grad{f_{\lambda}(x_k)}}^2
		+ \stepsize_{k-1}^2 \left(\nu L^2 + \thres^2/\lambda \right).		 
\end{align}
Finally, multiplying both sides of eq.\eqref{eq:lem:xdiff:non-smooth} by $1/(\lambda\nu)$ and adding the result to \eqref{eq:wcvx:lem:2:proof:envelope:3} (noting that $\rho \lambda \leq 1/2$) concludes the proof.

\section{Proof of Theorem~\ref{thrm:wcvx}}\label{appendix:proof:thrm:wcvx}

Taking the expectation on both sides of \eqref{eq:lem:Vfunc:non-smooth} and summing the result over $k=0,\ldots,K$ yields
\begin{align*}
		\frac{1}{2}	
		\sum_{k=0}^{K}	\stepsize_{k-1}
		\E \left[\ltwo{\grad{F_{\lambda}(x_{k})}}^2\right]
	 \leq 
		 {V_{-1}} - \E \left[V_{K}\right]		 
		+
		C \sum_{k=0}^{K}\stepsize_{k-1}^2.
\end{align*}
We can upper bound $V_{K}$ as
\begin{align*}
	V_{K} 
	&\geq
		f_\lambda(x_K) 
		- \frac{1}{2\nu}\ltwo{\grad{f_{\lambda}(x_{K})}}^2 
		+ \left(1+ \frac{1}{\lambda\nu}\right) f(x_K) 
	\nonumber\\
	&\geq
	  \left(2 + \frac{1}{\lambda\nu}\right) \inf_{x}f(x)
	  - \frac{2L^2}{\nu},
\end{align*}
where the first inequality holds since $f(x_K) \geq \inf_{x}f(x) = \inf_{x}f_{\lambda}(x)$, and the second one follows from the fact that $\ltwo{\grad{f_{\lambda}(x_{K})}} \leq 2L$.
As for $V_{-1}$, we have 
\begin{align*}
	V_{-1} 
	&= 
		f_\lambda(x_{-1}) 
		+ 
		\frac{1}{2\nu}\ltwo{d_{-1} - \grad{f_{\lambda}(x_{-1})}}^2 - \frac{1}{2\nu}\ltwo{\grad{f_{\lambda}(x_{-1})}}^2 
	\nonumber\\
		&\hspace{0.45cm}
		+ \left(1+\frac{1}{\lambda\nu} \right) f(x_{-1})
		+ \left(\frac{1-\mmt_{-1}}{2\lambda\nu^2} + \frac{\stepsize_{-1}}{\lambda\nu}\right)\ltwo{d_{-1}}^2
	\nonumber\\
	&\leq
		\left(2+\frac{1}{\lambda\nu} \right) f(x_{0})
\end{align*}
where we used the facts that $x_{-1}=x_0$,  $d_{-1}=0$, and $f_{\lambda}(x_{-1})=f_{\lambda}(x_{0}) \leq f(x_0)$. We thus arrive at
\begin{align*}
		\frac{1}{2}	
		\sum_{k=0}^{K}	\stepsize_{k-1}
		\E \left[\ltwo{\grad{F_{\lambda}(x_{k})}}^2\right]
	 \leq 
	 	\left(2+\frac{1}{\lambda\nu} \right) \Delta  + \frac{2L^2}{\nu}
		+
		C \sum_{k=0}^{K}\stepsize_{k-1}^2.
\end{align*}
Dividing both sides of the preceding inequality by $\sum_{k=0}^{K}	\stepsize_{k-1}=\sum_{i=0}^{K-1}	\stepsize_{i}$ (noting that $\stepsize_{-1}=0$) and using the definition of $k^*$ yields the first claim in the theorem.  
Finally, using the facts that $\nu=1/\stepsize_0$, $\beta_0 = 1/\sqrt{K}$, $\lambda=1/(2\rho)$, $\gamma\geq2L$, and  $\stepsize_0=1/\rho$, basic algebraic manipulations yield the last claim.

\end{document}